\documentclass[a4paper,11pt]{article}
\usepackage{amsmath,amsfonts,amssymb,amsthm}
\usepackage{enumitem}
\usepackage[stable]{footmisc}
\usepackage{graphicx}
\usepackage[all]{xy}
\usepackage{xcolor}
\usepackage[nonewpage]{imakeidx}

\makeindex

\def\spn{\bigskip\par\noindent}
\def\mpn{\medskip\par\noindent}
\def\pn{\par\noindent}

\def\SR{\mathsf{R}}
\def\SS{\mathsf{S}}
\def\Ss{\mathbf{S}}

\def\St{\mathbf{T}}

\def\CE{\mathcal{E}}
\def\CF{\mathcal{F}}
\def\CI{\mathcal{I}}
\def\CO{\mathcal{O}}
\def\CP{\mathcal{P}}

\def\CZ{\mathcal{Z}}
\def\Hom{\mathrm{Hom}}
\def\Mod{\hbox{-}\mathrm{Mod}}

\def\Id{\mathrm{Id}}
\def\Br{\mathrm{Br}}
\def\un{\mathbf{1}}
\def\Ind{\mathrm{Ind}}
\def\Inf{\mathrm{Inf}}

\def\Res{\mathrm{Res}}
\def\Iso{\mathrm{Iso}}
\def\Out{\mathrm{Out}}
\def\Aut{\mathrm{Aut}}
\def\Irr{\mathrm{Irr}}
\def\Inn{\mathrm{Inn}}

\def\F{\mathbb{F}}

\def\sfc{\mathsf{c}}
\def\normal{\trianglelefteq}
\definecolor{vert}{rgb}{0.12, 0.6, 0.17}
\newcommand{\Lu}[2]{#1\langle#2\rangle}

\def\sur{\overline}
\newcommand{\Sur}[1]{\varepsilon_\SR(#1)}

\def\Proj{{\mathsf{Proj}}}
\def\Cart{{\mathrm{Cart}}}
\def\smp{\smallskip\par}
\def\dom{\backslash}

\def\FF{\mathbb{F}}
\def\CC{\mathbb{C}}

\def\ZZ{\mathbb{Z}}
\def\NN{\mathbb{N}}
\newcommand{\pf}{{\flushleft\bf Proof: }}
\def\endpf{~\leaders\hbox to 1em{\hss\  \hss}\hfill~\raisebox{.5ex}{\framebox[1ex]{}}\smp}

\newcounter{nonce}[section]
\def\thenonce{\thesection.\arabic{nonce}}
\newenvironment{enonce}[1]{\pagebreak[2]\refstepcounter{nonce}\vskip 2ex plus 1ex minus 1ex\par\noindent{{\bf #1~\thesection.\arabic{nonce}:}}\begin{it}}{\end{it}\vskip 2ex plus 2ex minus 2ex\pn}
\newenvironment{enonce*}[1]{\pagebreak[2]\vskip 2ex plus 1ex minus 1ex\par\noindent{{\bf #1:}}\begin{it}}{\end{it}\vskip 2ex plus 2ex minus 2ex\pn}
\newenvironment{rem}[1]{\pagebreak[2]\refstepcounter{nonce}\vskip 2ex plus 1ex minus 1ex\par\noindent{{\bf #1~\thesection.\arabic{nonce}:}}}{\vskip 2ex plus 2ex minus 2ex\pn}
\newcounter{monequation}[nonce]
\newenvironment{moneq}{\refstepcounter{monequation}\begin{equation}}{\end{equation}}

\def\spn{\bigskip\par\noindent}
\def\pn{\par\noindent}
\def\mpn{\medskip\par\noindent}

\def\CE{\mathcal{E}}
\def\CF{\mathcal{F}}
\def\CP{\mathcal{P}}

\def\CT{\mathcal{T}}

\def\CZ{\mathcal{Z}}
\def\SS{\mathsf{S}}
\def\Hom{\mathrm{Hom}}
\def\Mod{\hbox{-}\mathrm{Mod}}

\def\Id{\mathrm{Id}}
\def\Br{\mathrm{Br}}
\def\un{\mathbf{1}}
\def\Ind{\mathrm{Ind}}
\def\Inf{\mathrm{Inf}}

\def\Res{\mathrm{Res}}
\def\Iso{\mathrm{Iso}}
\def\Out{\mathrm{Out}}
\def\Aut{\mathrm{Aut}}
\def\Inn{\mathrm{Inn}}

\def\F{\mathbb{F}}
\def\Z{\mathbb{Z}}
\def\normal{\trianglelefteq}
\definecolor{vert}{rgb}{0.12, 0.6, 0.17}
\def\sur{\overline}

\def\sfP{{\mathsf{P}}}
\def\Proj{{\mathsf{Proj}}}
\def\smp{\smallskip\par}
\def\dom{\backslash}

\def\FF{\mathbb{F}}
\def\CC{\mathbb{C}}
\def\NN{\mathbb{N}}

\def\Lu{L\langle u\rangle}
\def\uhat{{\langle u\rangle}^\natural}
\def\htheta{\hat{\theta}}
\def\dgd{{i_g}^\delta}
\def\semiG{\sur{G}_{Q,\delta}\ltimes\uhat}
\newcommand{\PimS}[1]{\mathrm{Pim}^\sharp(k\sur{N}_{#1})}
\newcommand{\Pim}[1]{\mathrm{Pim}(k\sur{N}_{#1})}
\newcommand{\ProjS}[1]{\F\mathrm{Proj}^\sharp(k\sur{N}_{#1})}
\def\findemo{~\leaders\hbox to 1em{\hss\  \hss}\hfill~\raisebox{.5ex}{\framebox[1ex]{}}\smp}
\newcommand{\tpar}{\medskip\par\noindent\pagebreak[3]\refstepcounter{nonce}{\bf \thenonce.\ }}
\newcommand{\ssur}[1]{\sur{#1}^{\,\flat}}

\newcommand{\ordidx}[1]{\makebox[0pt]{\phantom{#1}}}
\begin{document}
\centerline{\Large\bf Diagonal $p$-permutation functors in characteristic $p$}\vspace{1.5ex}\par
\centerline{\bf Serge Bouc and Deniz Y\i lmaz}\vspace{1.5ex}
\begin{abstract} Let $p$ be a prime number. We consider diagonal $p$-permutation functors over a (commutative, unital) ring $\SR$ in which all prime numbers different from $p$ are invertible. We first determine the finite groups~$G$ for which the associated essential algebra $\CE_\SR(G)$ is non zero: These are groups of the form $G=\Lu$, where $(L,u)$ is a $D^\Delta$-pair. \par
When $\SR$ is an algebraically closed field $\F$ of characteristic 0 or $p$, this yields a parametrization of the simple diagonal $p$-permutation functors over $\F$ by triples $(L,u,W)$, where $(L,u)$ is a $D^\Delta$-pair, and $W$ is a simple $\F\Out(L,u)$-module. \par
Finally, we describe the evaluations of the simple functor $\SS_{L,u,W}$ para\-metrized by the triple $(L,u,W)$. We show in particular that if $G$ is a finite group and $\F$ has characteristic $p$, the dimension of $\SS_{L,1,\F}(G)$ is equal to the number of conjugacy classes of $p$-regular elements of $G$ with defect isomorphic to $L$.
\end{abstract}

\section{Introduction} Let $k$ be an algebraically closed field of positive characteristic $p$, and $\SR$ be a commutative ring (with 1). As in~\cite{bouc-yilmaz} and~\cite{functorialequivalence}, we consider the following category category $\SR pp_k^\Delta$:
\begin{itemize}
\item The objects of $\SR pp_k^\Delta$ are the finite groups.
\item For finite groups $G$ and $H$, the hom set $\Hom_{\SR pp_k^\Delta}(G,H)$ is defined as $\SR T^\Delta(H,G)=\SR\otimes_\ZZ T^\Delta(H,G)$, where $T^\Delta(H,G)$ is the Grothendieck group of diagonal $p$-permutation $(kH,kG)$-bimodules.
\item The composition in $\SR pp_k^\Delta$ is induced by the usual tensor product of bimodules.
\item The identity morphism of the group $G$ is the class of the $(kG,kG)$-bimodule $kG$.
\end{itemize}
The category $\CF_{\SR pp_k}^\Delta$ of {\em diagonal $p$-permutation functors over $\SR$} is the category of $\SR$-linear functors from $\SR pp_k^\Delta$ to the category $\SR\Mod$ of all $\SR$-modules. It is an abelian category.\par
In~\cite{bouc-yilmaz} and~\cite{functorialequivalence}, we mainly considered the case where $\SR$ is an algebraically closed field $\FF$ of characteristic 0. In~\cite{functorialequivalence}, we showed in particular that the category $\CF_{\FF pp_k}^\Delta$ is semisimple, and we classified and described its simple objects. For an arbitrary commutative ring $\SR$, we also introduced a new equivalence for blocks of groups algebras, called {\em functorial equivalence} over~$\SR$, using diagonal $p$-permutation functors over $\SR$ naturally attached to pairs $(G,b)$ of a finite group $G$ and a block idempotent $b$ of $kG$.\par
A natural question is then to see what happens when $\SR$ is a field of {\em positive characteristic}, in particular when $\SR=k$. The main result of the present paper is a parametrization and a description of the simple diagonal $p$-permutation functors in characteristic $p$. \par
A key ingredient to the parametrization and description of simple functors is {\em the essential algebra} $\CE_\SR(G)$ of a group $G$, namely the quotient of the endomorphism algebra $\SR T^\Delta(G,G)$ of $G$ in the category $\SR pp_k^\Delta$ by the ideal of linear combinations of endomorphisms which factor through a group of order strictly smaller than $|G|$. We first find some conditions on $G$ (Corollary~\ref{semidirect}, Theorem~\ref{K cyclic}, Theorem~\ref{K faithful}) for the {(non-)vanishing} of $\CE_\SR(G)$. In particular, if any prime number different from $p$ is invertible in $\SR$, we show (Corollary~\ref{K cyclic faithful}) that $\CE_\SR(G)$ is non zero if and only if $G$ is a semidirect product $L\rtimes \langle u\rangle$ (which we denote by $\Lu$), where $(L,u)$ is a {\em $D^\Delta$-pair}, that is a pair of a finite $p$-group $L$ and a $p'$-automorphism $u$ of $L$ (Definition~\ref{D Delta pair}). Moreover, we describe completely (Theorem~\ref{structure essential}) the structure of the algebra $\CE_\SR(G)$ in this case.\par
In Section~5, we study simple diagonal $p$-permutation functors, so we assume that $\SR$ is a field $\F$ of characteristic 0 or $p$. Applying the results of the previous sections, we know that if $S$ is a simple diagonal $p$-permutation functor over $\F$, then a minimal group for $S$ is of the form $\Lu$, where $(L,u)$ is a $D^\Delta$-pair, and the evaluation $V=S\big(\Lu\big)$ is a simple $\CE_\F(\Lu)$-module. Conversely, to any triple $(L,u,V)$, where $(L,u)$ is a $D^\Delta$-pair and $V$ is a simple $\CE_\F(\Lu)$-module, we associate a simple functor $S_{\Lu,V}$ with minimal group $\Lu$, and such that $S_{\Lu,V}\big(\Lu\big)\cong V$. Then we compute (Theorem~\ref{isomorphism}) the evaluation $S_{\Lu,V}(G)$ at an arbitrary finite group $G$.\par
The precise structure of the essential algebra given by Theorem~\ref{structure essential} now allows for another parametrization of the simple functors, namely by triples $(L,u,W)$, where $(L,u)$ is a $D^\Delta$-pair and $W$ is a simple module for the algebra $\F\Out(L,u)$-module of the group $\Out(L,u)$  (Notation~\ref{Out(L,u)}) of outer automorphisms of $(L,u)$. In Theorem~\ref{the simple functors}, we describe the evaluations of the simple functor $\SS_{L,u,W}$ parametrized by such a triple $(L,u,W)$.\par
Section 6 is devoted to some examples: First the simple functor $S_{\un,1,\F}$, which turns out to be closely related to the Cartan map (Lemma~\ref{unique minimal}, Proposition~\ref{basis}). This example shows in particular that the category $\CF_{\F pp_k}$ is {\em not} semisimple when $\F$ has characteristic $p$. Then we describe (Theorem~\ref{simple L,1,W}) the evaluations of the simple functor $\SS_{L,1,W}$. In particular (Corollary~\ref{simple L,1,k}), we show that for a finite group $G$, the dimension of $\SS_{L,1,k}(G)$ is equal to the number \nopagebreak of conjugacy classes of $p'$-elements of $G$ with {\em defect isomorphic to $L$} (Definition~\ref{defect isomorphic}).
\section{Notation and terminology\footnote{An additional list of symbols is included at the end of the paper.}}
Throughout the paper:
\begin{itemize}[label=$\blacktriangleright$,leftmargin=3ex,itemsep=0.5ex]
\item $k$ is an algebraically closed field of positive characteristic $p$.
\item $\SR$ is a commutative ring (with 1).
\item For a finite group $G$, we denote by $\Proj(kG)$ the group of projective $kG$-modules, and by $R_k(G)$ the Grothendieck group of the category of finite dimensional $kG$-modules. We set $\SR \Proj(G)=\SR\otimes_\Z\Proj(G)$ and $\SR R_k(G)=\SR\otimes_\Z R_k(G)$.
\item If $P$ is a $p$-subgroup of a finite group $G$, and $M$ is a $kG$-module, we denote by $M[P]$ the Brauer quotient of $M$ at $P$, and by $\Br_P:M^P\to M[P]$ the projection map. The module $M[P]$ is a $k\sur{N}_G(P)$-module, where $\sur{N}_G(P)=N_G(P)/P$.
\item For a finite group $G$, a {\em $p$-permutation $kG$-module} (see~\cite{brouepperm}) is a direct summand of a permutation $kG$-module, i.e. of a module admitting a $G$-invariant $k$-basis. Equivalently, a $kG$-module $M$ is a $p$-permutation module if the restriction $\Res_S^GM$ of $M$ to a Sylow $p$-subgroup $S$ of $G$ is a permutation $kS$-module.
\item From~\cite{brouepperm}, we know that the indecomposable $p$-permutation $kG$-modules (up to isomorphism) are parametrized by pairs $(P,E)$, where $P$ is a $p$-subgroup of $G$, up to conjugation, and $E$ is an indecomposable projective $k\sur{N}_G(P)$-module, up to isomorphism. The indecomposable module $M(P,E)$ para\-metrized by the pair $(P,E)$ is the only indecomposable direct summand  with vertex $P$ of $L_{P,E}=\Ind_{N_G(P)}^G\Inf_{\sur{N}_G(P)}^{N_G(P)}E$, the other direct indecomposable summands having vertex strictly contained in $P$, up to conjugation. The module $M(P,E)$ has vertex $P$, and $M(P,E)[P]\cong E$ as $k\sur{N}_G(P)$-modules.
\item It follows that the Grothendieck group of $p$-permutation $kG$-modules, for relations given by direct sum decomposition, has a basis consisting of the modules $L_{P,E}$, where $P$ is a $p$-subgroup of $G$, up to conjugation, and $E$ is an indecomposable projective $k\sur{N}_G(P)$-module.
\item When $G$ and $H$ are finite groups, and $L$ is a subgroup of $H\times G$, we denote by $p_1(L)$ (resp. $p_2(L)$) the projection of $L$ in $H$ (resp. in $G$), and we set
$$k_1(L)=\{h\in H\mid (h,1)\in L\}\hspace{2ex}\hbox{and}\hspace{2ex}k_2(L)=\{g\in G\mid (1,g)\in L\}.$$
We say that $L$ is {\em diagonal} if $k_1(L)=k_2(L)=1$. Equivalently,
$$L=\Delta(Y,\pi,X)=\big\{\big(\pi(x),x\big)\mid x\in X\big\},$$
where $X$ is a subgroup of $G$ and $\pi:X\to Y$ is a group isomorphism from $X$ to a subgroup $Y$ of $H$. If $X=Y$ and $\pi=\Id$, we simply write $\Delta(X)=\Delta(X,\Id,X)$. For $X\leq G$ and an embedding $\psi:X\hookrightarrow H$, we also write $\Delta_\psi(X)$ instead of $\Delta\big(\psi(X),\psi,X\big)$.
\item For finite groups $G$ and $H$, a $p$-permutation $(kH,kG)$-bimodule is a $(kH,kG)$-bimodule which is a $p$-permutation module when viewed as a $k(H\times G)$-module. A $p$-permutation $(kH,kG)$-bimodule $M$ is {\em diagonal} if in addition $M$ is projective when viewed as a left $kH$-module and a right $kG$-module. Equivalently $M$ is a $p$-permutation $(kH,kG)$-bimodule, and all the vertices of the indecomposable summands of $M$ are diagonal $p$-subgroups of $H\times G$. 
\item For finite groups $G$ and $H$, we denote by $T^\Delta(H,G)$ the Grothendieck group of diagonal $p$-permutation $(kH,kG)$-bimodules, for relations given by direct sum decomposition. We set $\SR T^\Delta(H,G)=\SR\otimes_\ZZ T^\Delta(H,G)$. The group $T^\Delta(H,G)$ has a basis consisting of the bimodules of the form 
$$\Ind_{N_{H\times G}(P)}^{H\times G}\Inf_{\sur{N}_{H\times G}(P)}^{N_{H\times G}(P)}E,$$
where $P$ is a diagonal $p$-subgroup of $H\times G$ (up to conjugation), and $E$ is an indecomposable projective $\sur{N}_{H\times G}(P)$-module.
\item When $G$, $H$, and $K$ are finite groups, if $M$ is a diagonal $p$-permutation $(kG,kH)$-bimodule and $N$ is a diagonal $p$-permutation $(kK,kH)$-bimo\-dule, then $N\otimes_{kH}M$ is a diagonal $p$-permutation $(kK,kG)$-bimodule. This induces a well defined bilinear map 
$$T^\Delta(K,H)\times T^\Delta(H,G)\to T^\Delta(K,G),$$
still denoted $(v,u)\mapsto v\otimes_{kH}u$. This bilinear map is also the composition in the category $\SR pp_k^\Delta$ of the introduction, so it will be sometimes denoted by $(v,u)\mapsto v\circ u$.
\item For finite groups $G$ and $H$, we say that an element $u\in \SR T^\Delta(G,H)$ is {\em right essential} (resp. {\em left essential}) if it cannot be factored through groups of order strictly smaller than $|H|$ (resp. of order strictly smaller than $|G|$), that is if $u\notin\sum_{|K|<|H|}\limits\SR T^\Delta(G,K)\circ \SR T^\Delta(K,H)$ (resp. if $u\notin\sum_{|K|<|G|}\limits\SR T^\Delta(G,K)\circ \SR T^\Delta(K,H)$). A $(kG,kH)$-bimodule $M$ is called right essential over $\SR$ - or simply right essential - (resp. left essential) if the element $M$ of $\SR T^\Delta(G,H)$ is right essential (resp. left essential). \par
If $|G|=|H|$, being left essential is equivalent to being right essential, so we simply say {\em essential}.
\item In particular, for a finite group $G$, the endomorphism algebra of $G$ in the category $\SR pp_k^\Delta$ is $\SR T^\Delta(G,G)$. The {\em essential algebra} (over $\SR$) of $G$ is the quotient
$$\CE_\SR(G)\index{\ordidx{EA}$\CE_\SR(G)$}=\SR T^\Delta(G,G)/\sum_{|H|<|G|}\SR T^\Delta(G,H)\circ \SR T^\Delta(H,G)$$
of $\SR T^\Delta(G,G)$ by the (two sided) ideal of non-essential elements. We denote by $u\mapsto \Sur{u}$ the projection map $\SR T^\Delta(G,G)\to \CE_\SR(G)$.
\item The main reason for considering the previous essential algebra is the following: By standard results (see e.g.~\cite{corfun-finiteness}, Lemma~2.5 and Proposition~2.7), if $S$ is a simple diagonal $p$-permutation functor over $\SR$, and if $G$ is a group such that $V:=S(G)\neq 0$, then $V$ is a simple $\SR T^\Delta(G,G)$-module, and $S$ is isomorphic to the unique simple quotient $S_{G,V}$ of the functor $L_{G,V}:H\mapsto \SR T^\Delta(H,G)\otimes_{\SR T^\Delta(G,G)}V$. Moreover, if $G$ is a group of minimal order such that $S(G)\neq 0$, then in fact $V$ is a simple $\CE_\SR(G)$-module, and $S\cong S_{G,V}$. So we are looking for pairs $(G,V)$ of a finite group $G$ and a simple $\CE_\SR(G)$-module. In particular, for such a pair, the essential algebra $\CE_\SR(G)$ is non-zero.
\item An elementary group (or Brauer elementary group) is a finite group of the form $Q\times C$, where $Q$ is a $q$-group for some prime number $q$, and $C$ is a cyclic group (that can be assumed of order prime to $q$). When $p$ is a prime number, an elementary $p'$-group is an elementary group of order prime to~$p$.
\end{itemize}
\section{Vanishing of $\CE_\SR(G)$}
Let $G$ be a finite group. We want to know when the essential algebra $\CE_\SR(G)$ is non-zero. We start with some classical lemmas.
\begin{enonce}{Lemma} \label{inner}Let $G=P\rtimes K$, where $P$ and $K$ are finite groups of coprime order. Let moreover $\varphi$ be an automorphism of $G$. Then:
\begin{enumerate}
\item $\varphi(P)=P$.
\item $C_G(P)=Z(P)C_K(P)$. 
\item Suppose that $C_G(P)=Z(P)$, or equivalently by Assertion 2, that $K$ acts faithfully on $P$. Then the following are equivalent:
\begin{enumerate}
\item The restriction of $\varphi$ to $P$ is the identity.
\item $\varphi$ is an inner automorphism $i_w$ of $G$, for some $w\in Z(P)$.
\end{enumerate}
\end{enumerate}
\end{enonce}
\pf 1. This is clear, since $P$ is the set of elements of $G$ of order dividing the order of $P$.\mpn
2. The inclusion $Z(P)C_K(P)\leq C_G(P)$ is clear. Conversely, if $xt\in C_G(P)$, where $x\in P$ and $t\in K$, then $^ty=y^x$ for any $y\in P$. It follows that $^{t^n}y=y^{x^n}$ for any $n\in \NN$ and any $y\in P$. Taking $n=|x|$ gives $t^n\in C_K(P)$, hence $t\in C_K(P)$ since $(n,|t|)=1$. Then $x\in C_P(P)=Z(P)$.\mpn
3. It is clear that $(b)$ implies $(a)$\footnote{and we don't need the assumption $C_G(P)=Z(P)$ for that\ldots}. For the converse, assume that $(a)$ holds, and that $C_G(P)=Z(P)$. Let $x,y\in P$ and $s,t\in K$. Then
\begin{align*}
\varphi(xs\cdot yt)&=\varphi(x\,{^sy}\cdot st)=x\,{^sy}\varphi(st)\\
&=\varphi(xs)\varphi(yt)=x\varphi(s)y\varphi(t)=x\,{^{\varphi(s)}y}\,\varphi(s)\varphi(t).
\end{align*}
Hence $^sy={^{\varphi(s)}y}$ for any $y\in P$. In other words $z(s):=s^{-1}\varphi(s)\in C_G(P)$, so $z$ is a map from $K$ to $Z(P)$. \par
Now $\varphi(s)=sz(s)$, so $z(st)=z(s)^t\,z(t)$ for any $s,t\in K$. In other words, the map $z$ is a crossed morphism from $K$ to $Z(P)$. Since $K$ and $Z(P)$ have coprime order, it follows that there exists $w\in Z(P)$ such that $z(s)=w^s\cdot w^{-1}$, for any $s\in K$. In other words $\varphi(s)=s\cdot w^s\cdot w^{-1}=wsw^{-1}=i_w(s)$. Since $i_w(x)=x=\varphi(x)$ for any $x\in P$, it follows that $\varphi=i_w$.\findemo 

\begin{enonce}{Lemma} \label{kG bimodule}Let $G$ be a finite group, and $P$ be a normal $p$-subgroup of $G$. Then:
\begin{enumerate}
\item $N_{G\times G}\big(\Delta(P)\big)=\{(a,b)\in G\times G\mid ab^{-1}\in C_G(P)\}$.
\item  Set $N=N_{G\times G}\big(\Delta(P)\big)$ and $\sur{N}=N/\Delta(P)$. There is an isomorphism of $(kG,kG)$-bimodules
$$kG\cong\Ind_{N}^{G\times G}\Inf_{\sur{N}}^NkC_G(P),$$
where the action of $\sur{N}$ on $kC_G(P)$ is given by $(a,b)\Delta(P)\cdot\gamma=a\gamma b^{-1}$.
\end{enumerate}
\end{enonce}
\pf 1. This is clear, since $(a,b)\in N$ if and only if $x^a=x^b$, i.e. $x^{ab^{-1}}=x$, for all $x\in P$.\spn
2. The group $\sur{N}$ permutes the set $C_G(P)$ transitively, and the stabilizer in $\sur{N}$ of $1\in C_G(P)$ is the group $\{(a,a)\Delta(P)\mid a\in G\}=\Delta(G)/\Delta(P)$. So $kC_G(P)\cong \Ind_{\Delta(G)/\Delta(P)}^{\sur{N}}k$, and
\begin{align*}
\Ind_{N}^{G\times G}\Inf_{\sur{N}}^NkC_G(P)&\cong \Ind_{N}^{G\times G}\Inf_{\sur{N}}^N\Ind_{\Delta(G)/\Delta(P)}^{\sur{N}}k\\
&\cong \Ind_{N}^{G\times G}\Ind_{\Delta(G)}^N\Inf_{\Delta(G)/\Delta(P)}^{\Delta(G)}k\\
&\cong \Ind_{\Delta(G)}^{G\times G}k\cong kG
\end{align*}
as $(kG,kG)$-bimodules.\endpf
The next step is an important reduction allowed by the following stronger version of Dress induction theorem, due to Boltje and K\"ulshammer (\cite{boltje-kulshammer-dressinduction},~Theorem 3.3): 
\begin{enonce}{Theorem} \label{BK}Let $H$ be a finite group, and $U$ be an indecomposable $kH$-module with vertex $D$ and source $Z$. Then, in the Green ring of $kH$, we have
$$[U]=\sum_{i=1}^na_i[\Ind_{H_i}^GV_i],$$
where, for $i=1,\ldots,n$:
\begin{itemize}
\item $a_i$ is an integer.
\item $H_i$ is a subgroup of $H$ such that $D_i:=O_p(H_i)\leq D$ and $H_i/D_i$ is an elementary $p'$-group.
\item $V_i$ is an indecomposable $kH_i$-module with vertex $D_i$ and source $\Res_{D_i}^{H_i}V_i$, which is a direct summand of $\Res_{D_i}^D Z$.
\end{itemize}
\end{enonce}
\begin{enonce}{Corollary} \label{semidirect}Let $G$ be a finite group. Then $\CE_\SR(G)=0$ unless $G\cong P\rtimes K$, where $P$ is a $p$-group, and $K$ is an elementary $p'$-group.
\end{enonce}
\pf We apply Theorem~\ref{BK} to the case $H=G\times G$ and $U=kGb$, where $b$ is a block idempotent of $kG$. Then $U$ is a diagonal $p$-permutation $(kG,kG)$-bimodule with diagonal vertex $\Delta(D)=\Delta(D,\Id,D)$, where $D\leq G$ is a defect group of $b$. We can conclude that $[U]$ is a linear combination with integer coefficients of (isomorphism classes of) induced bimodules $\big[\Ind_{H_i}^{G\times G}V_i\big]$, where $H_i$ is a subgroup of $G\times G$ such that $D_i=O_p(H_i)\leq \Delta(D)$ and $H_i/O_p(H_i)$ is an elementary $p'$-group. Let $G_i$ be the first projection of $H_i$ on $G$. Then the bimodule $\Ind_{H_i}^{G\times G}V_i$ factors as
$$\Ind_{H_i}^{G\times G}V_i\cong \Ind_{\Delta(G_i)}^{G\times G_i}\otimes_{kG_i}\Ind_{H_i}^{G_i\times G}V_i,$$
where $H_i$ on the right hand side is viewed as a subgroup of $G_i\times G$. Now if $G_i<G$, then the image of $kGb$ in $\CE_\SR(G)$ is equal to 0. And if $G_i=G$, then $G$ is a quotient of $H_i$, so $G/O_p(G)$ is an elementary $p'$-group. In other words $G\cong P\rtimes K$, where $P$ is a $p$-group, and $K$ is an elementary $p'$-group. \par
Now $\CE_\SR(G)$ is non zero if and only if its identity element is non zero, that is if the image of the bimodule $kG$ in $\CE_\SR(G)$ is non zero. Since $kG$ is the direct sum of the bimodules $kGb$, when $b$ runs through block idempotents of $kG$, there is at least one such idempotent $b$ such that the image of $kGb$ in $\CE_\SR(G)$ is non-zero. Hence $G\cong P\rtimes K$, where $P$ is a $p$-group and $K$ is an elementary $p'$-group.\endpf
\begin{enonce}{Lemma} \label{does not factor}Let $G=P\rtimes K$, where $P$ is a finite $p$-group, and $K$ is an elementary $p'$-group. Let $H$ be a finite group, and $U$ be a right essential indecomposable diagonal $p$-permutation $(kG,kH)$-bimodule. Then:
\begin{enumerate}
\item The essential algebra $\CE_\SR(H)$ is non-zero. In particular $H=Q\rtimes L$, where $Q$ is a $p$-group and $L$ is an elementary $p'$-group.
\item There exist an injective group homomorphism $\pi:Q\hookrightarrow P$, a subgroup $T$ of $N_{K\times L}\big(\Delta_\pi(Q)\big)$ with $p_2(T)=L$, and a simple $kT$-module $W$ such that
$$U\cong\Ind_{\Delta_\pi(Q)\cdot T}^{G\times H}\Inf_{T}^{\Delta_\pi(Q)\cdot T}W$$
as $(kG,kH)$-bimodules.
\end{enumerate}
\end{enonce}
\pf  1. If $\CE_\SR(H)=0$, the identity $(kH,kH)$-bimodule $kH$ factors through groups of order stritcly smaller than $|H|$, so the same holds for $U$ (by right composition with $kH$). Hence $\CE_\SR(H)\neq 0$, and Assertion 1 follows from~Corollary~\ref{semidirect}. \mpn
2. From Assertion 1 follows in particular that the group $G\times H$ is solvable, with a normal Sylow $p$-subgroup $P\times Q$. Then there is a diagonal $p$-subgroup $\Delta_\pi(S)$ of $G\times H$, where $S$ is a $p$-subgroup of $H$ (that is, a subgroup of $Q$) and $\pi:S\hookrightarrow P$ is an injective group homomorphism, such that
$$U\cong\Ind_N^{G\times H}\Inf_{\sur{N}}^NE,$$
where $N=N_{G\times H}\big(\Delta_\pi(S)\big)$ and $\sur{N}=N/\Delta_\pi(S)$, and $E$ is an indecomposable projective $k\sur{N}$-module. \par
Now $\sur{N}$ itself also has a normal Sylow $p$-subgroup $X$, and there is a $p'$-subgroup $\sur{T}$ of $\sur{N}$ such that $\sur{N}=X\rtimes \sur{T}$. Moreover $\sur{T}$ lifts to a $p'$-subgroup $T$ of $N$, that we can assume contained in the $p'$-Hall subgroup $K\times L$ of $G\times H$, up to replacing $\Delta_\pi(S)$ by a conjugate subgroup. Finally $E\cong\Ind_{\sur{T}}^{\sur{N}}\sur{W}$, where $\sur{W}$ is a simple $k\sur{T}$-module. Let $W$ be the simple $kT$-module corresponding to $\sur{W}$ via the isomorphism $\sur{T}\cong T$. It follows that
\begin{align*}
U&\cong \Ind_N^{G\times H}\Inf_{\sur{N}}^N\Ind_{\sur{T}}^{\sur{N}}\sur{W}\\
&\cong \Ind_N^{G\times H}\Ind_{\Delta_\pi(S)\cdot T}^N\Inf_{T}^{\Delta_\pi(S)\cdot T}W\\
&\cong\Ind_{\Delta_\pi(S)\cdot T}^{G\times H}\Inf_{T}^{\Delta_\pi(S)\cdot T}W.
\end{align*}
Since $U$ is right essential, we have $p_2\big(\Delta_\pi(S)\cdot T\big)=H=Q\cdot L$. This forces $S=Q$, and $p_2(T)=L$, proving Assertion 2.\endpf
\begin{enonce}{Theorem} \label{K cyclic} Let $G$ be a finite group of the form $G=P\rtimes K$, where $P$ is a $p$-group and $K$ is a non-cyclic elementary $p'$-group. Then $|K|^2\,\CE_\SR(G)=0$. In particular, if $|K|$ is invertible in $\SR$, then $\CE_\SR(G)=0$.
\end{enonce}
\pf Let $M$ be an essential indecomposable diagonal $p$-permutation $(kG,kG)$-bimodule. By Lemma~\ref{does not factor}, applied to $H=G$ and $U=M$, we know that 
$$M\cong \Ind_{\Delta_\pi(P)\cdot T}^{G\times G}\Inf_{T}^{\Delta_\varphi(P)\cdot T}W,$$
for some $\pi\in \Aut(P)$, some subgroup $T$ of $N_{K\times K}\big(\Delta_\pi(P)\big)$ with $p_2(T)=K$, and some simple $kT$-module $W$. \par
Now $T\leq K\times K$, so $T$ is a $p'$-group, and $|T|$ divides $|K|^2$. Moreover by Artin's induction theorem, in $R_k(T)\cong R_\CC(T)$, we have an equality of the form
$$|T|W=\sum_{i=1}^nn_i\Ind_{C_i}^Tk_{\lambda_i},$$
where, for $1\leq i\leq n$, $C_i$ is a cyclic subgroup of $T$, $n_i$ is an integer, and $k_{\lambda_i}$ is a one dimensional $kC_i$-module. Hence in $\SR T^\Delta(G,G)$, we have that
\begin{align*}
|T|M&=\sum_{i=1}^nn_i\Ind_{\Delta_\pi(P)\cdot T}^{G\times G}\Inf_T^{\Delta_\pi(P)\cdot T}\Ind_{C_i}^Tk_{\lambda_i}\\
&=\sum_{i=1}^nn_i\Ind_{\Delta_\pi(P)\cdot T}^{G\times G}\Ind_{\Delta_\pi(P)\cdot C_i}^{\Delta_\pi(P)\cdot T}\Inf_{C_i}^{\Delta_\pi(P)\cdot C_i}k_{\lambda_i}\\
&=\sum_{i=1}^nn_i\Ind^{G\times G}_{\Delta_\pi(P)\cdot C_i}\Inf_{C_i}^{\Delta_\pi(P)\cdot C_i}k_{\lambda_i}.\\
\end{align*}
The image of $|K|^2\,M$ in $\CE_\SR(G)$ is equal to the image of
$$\frac{|K|^2}{|T|}\sum_{i=1}^nn_i\Ind^{G\times G}_{\Delta_\pi(P)\cdot C_i}\Inf_{C_i}^{\Delta_\pi(P)\cdot C_i}k_{\lambda_i},$$
which is equal to zero unless there exists $i\in\{1,\ldots,n\}$ such that 
$$p_1\big(\Delta_\varphi(P)\cdot C_i\big)=p_2\big(\Delta_\varphi(P)\cdot C_i\big)=G=P\cdot K.$$
This implies that $K$ is a quotient of $C_i$. Hence $K$ is cyclic, which completes the proof.\endpf
\begin{enonce}{Theorem} \label{K faithful}Let $G=P\rtimes K$, where $P$ is a $p$-group and $K$ is a cyclic $p'$-group. If $C_K(P)\neq 1$, then $\CE_\SR(G)=0$.
\end{enonce}
\pf\footnote{This proof is a slightly simplified and generalized version of the proof given by M.~Ducellier in Proposition~4.1.2 of his thesis~\cite{ducellier} in the case $\SR=\CC$.} Since $G$ has a normal Sylow $p$-subgroup, all the blocks of $G$ have defect $P$. Moreover, if $b$ is a block idempotent of $G$, then $b$ is a linear combination of $p$-regular elements of $C_G(P)=Z(P)\times C_K(P)$, so $b\in kC_K(P)$. \par
Since $C_K(P)\leq Z(G)$, it means that the block idempotents of $kG$ are exactly the primitive idempotents of the (split semisimple commutative) algebra $kC_K(P)$. Let $e$ be one of them, and $k_\lambda=kC_K(P)e$ be the corresponding (one dimensional) simple $kC_K(P)$-module, where $\lambda:C_K(P)\to k^\times$ is the associated group homomorphism.\par
Let $\varpi:G\to K$ denote the projection map. Set $N=N_{G\times G}\big(\Delta(P)\big)$ and $\sur{N}=N/\Delta(P)$. There is a short exact sequence
$$\xymatrix{1\ar[r]& Z(P)\ar[r]^-i&\sur{N}\ar[r]^-s&\widetilde{K}\ar[r]&1,}$$
where \begin{itemize}
\item $\widetilde{K}=\{(a,b)\in K\times K\mid a^{-1}b\in C_K(P)\}$. 
\item $i$ is the map sending $z\in Z(P)$ to $(z,1)\Delta(P)\in\sur{N}$. 
\item $s$ is the map sending $(a,b)\Delta(P)$ to $\big(\varpi(a),\varpi(b)\big)$.
\end{itemize}
So $\sur{N}\cong Z(P)\rtimes \widetilde{K}$, with the explicit embedding $\widetilde{K}\to \sur{N}$ sending $(a,b)\in \widetilde{K}$ to $(a,b)\Delta(P)\in\sur{N}$. We consider $\widetilde{K}$ as a subgroup of $\sur{N}$ via this embedding.\medskip\par
The Brauer quotient of the $(kG,kG)$-bimodule $kG$ at $\Delta(P)$ is isomorphic to $kC_G(P)$, so $kGe[\Delta(P)]\cong kC_G(P)\Br_P(e)=kC_G(P)e=kZ(P)\otimes_k k_\lambda$, since $\Br_P(e)=e$ as $e\in kC_K(P)$. It follows that 
$$kGe\cong \Ind_N^{G\times G}\Inf_{\sur{N}}^NkC_G(P)e\cong \Ind_N^{G\times G}\Inf_{\sur{N}}^N \big(kZ(P)\otimes_k k_\lambda\big),$$
where $\sur{N}\cong Z(P)\rtimes \widetilde{K}$ acts on $kZ(P)\otimes_k k_\lambda$ by
$$(a,b)\Delta(P).(z\otimes 1)=(ab^{-1})_p\,z\otimes (ab^{-1})_{p'}\cdot1=\lambda\big((ab^{-1})_{p'}\big)(ab^{-1})_p\,z\otimes 1$$
for $(a,b)\in N$ and $z\in Z(P)$, where $(ab^{-1})_p\in Z(P)$ and $(ab^{-1})_{p'}\in C_K(P)$ are the $p$-part and $p'$-part of $ab^{-1}\in C_G(P)=Z(P)\times C_K(P)$, respectively. Then $kZ(P)\otimes _kk_\lambda$ is isomorphic to $\Ind_{\widetilde{K}}^{\sur{N}}\,k_{\widetilde{\lambda}}$, where $\widetilde{\lambda}:\widetilde{K}\to k^\times$ sends $(a,b)\in \widetilde{K}$ to $\lambda(ab^{-1})\in k^\times$. So
\begin{align}\refstepcounter{monequation}
kGe&\cong\Ind_{N}^{G\times G}\Inf_{\sur{N}}^{N}\Ind_{\widetilde{K}}^{\sur{N}}k_{\widetilde{\lambda}}\nonumber\\
&\cong \Ind_{N}^{G\times G}\Ind_{\Delta(P)\widetilde{K}}^{N}\Inf_{\widetilde{K}}^{\Delta(P)\widetilde{K}}k_{\widetilde{\lambda}}\nonumber\\
&\cong \Ind^{G\times G}_{\Delta(P)\widetilde{K}}\Inf_{\widetilde{K}}^{\Delta(P)\widetilde{K}}k_{\widetilde{\lambda}}\label{kGe}.
\end{align}\par
Now we set $\sur{G}=G/C_K(P)$. We denote by $g\mapsto \sur{g}$ the projection map, and by $\delta:G\to G\times\sur{G}$ the map $g\mapsto (g,\sur{g})$. We will show that the $(kG,kG)$-bimodule $kGe$ factors through the group~$\sur{G}$, that is, there is a diagonal $p$-permutation $(kG,k\sur{G})$-bimodule $U$ and a diagonal $p$-permutation $(k\sur{G},kG)$-bimodule $V$ such that $kGe\cong U\otimes_{k\sur{G}}V$.\smallskip\par
The group $P$ embeds in $G\times\sur{G}$ via $\delta$. Its image $\delta(P)$ is  a diagonal subgroup of $G\times \sur{G}$, and its normalizer is
$$N_\delta:=N_{G\times\sur{G}}\big(\delta(P)\big)=\big\{(a,\sur{b})\in G\times \sur{G}\mid \sur{x^a}=\sur{x^b},\;\forall x\in P\big\}.$$
In other words $(a,\sur{b})\in N_\delta$ if and only if $\sur{x^{ab^{-1}}}=\sur{x}$ for all $x\in P$, or equivalently if the commutator $[x,ab^{-1}]$ is in $C_K(P)$. But since $P\normal G$, we have that $[P,ab^{-1}]\subseteq P$. Hence $(a,\sur{b})$ normalizes $\delta(P)$ if and only if $[P,ab^{-1}]\subseteq P\cap C_K(P)=1$, i.e. $ab^{-1}\in C_G(P)$. Thus
$$N_\delta=\{(a,\sur{b})\in G\times \sur{G}\mid ab^{-1}\in Z(P)\times C_K(P)\}.$$
Recall that $\varpi:G\to K$ denotes the projection map. We have a surjective group homomorphism $\sigma:N_\delta\to K$ sending $(a,\sur{b})$ to $\varpi(a)$. It induces a surjective group homomorphism
$$\sur{\sigma}:\sur{N}_\delta:=N_\delta/\delta(P)\to K$$
sending $(a,\sur{b})\delta(P)$ to $\varpi(a)$. The kernel of this morphism consists of the elements $(a,\sur{b})\delta(P)$ such that $a\in P$ and $ab^{-1}\in Z(P)\times C_K(P)$. Since 
$$(a,\sur{b})\delta(P)=(1,\sur{ba^{-1}})(a,\sur{a})\delta(P)=(1,\sur{ba^{-1}})\delta(P),$$
and since $\sur{ba^{-1}}\in C_G(P)/C_K(P)=Z(P)$, we get a short exact sequence
$$\xymatrix{
1\ar[r]&Z(P)\ar[r]^-\iota&\sur{N}_\delta\ar[r]^-{\sur{\sigma}}&K\ar[r]&1
}
$$
where $\iota(z)=(1,\sur{z})\delta(P)$ for $z\in Z(P)$. This sequence is split, via the morphism $a\in K\mapsto (a,\sur{a})\delta(P)$, for $a\in K$, so $\sur{N}_\delta\cong Z(P)\rtimes K$.\par
Now since $K$ is cyclic, we can extend $\lambda:C_K(P)\to k^\times$ to a group homomorphism $\beta:K\to k^\times$. This gives a one dimensional $kK$-module $k_\beta$, that we can induce to $\sur{N}_\delta=Z(P)\rtimes K$. We get a projective $k\sur{N}_\delta$-module, and we set
$$U:=\Ind_{N_\delta}^{G\times \sur{G}}\Inf_{\sur{N}_\delta}^{N_\delta}\Ind_{K}^{\sur{N}_\delta}k_\beta.$$
This is a diagonal $p$-permutation $(kG,k\sur{G})$-bimodule, and
$$U\cong \Ind_{N_\delta}^{G\times \sur{G}}\Ind_{\delta(P)K}^{N_\delta}\Inf_{K}^{\delta(P)K}k_\beta\cong \Ind_{\delta(P)K}^{G\times\sur{G}}\Inf_{K}^{\delta(P)K}k_\beta,$$
where $K$ is viewed as a subgroup of $N_\delta$ via the map $a\in K\mapsto (a,\sur{a})\in N_\delta$. We observe that $\delta(P)K$ is equal to $\delta(G)$, so
$$U\cong \Ind_{\delta(G)}^{G\times\sur{G}}\Inf_{K}^{\delta(G)}k_\beta.$$
We define similarly a $(k\sur{G},kG)$-bimodule $V$ by
$$V:=\Ind_{\delta^o(G)}^{\sur{G}\times G}\Inf_{K^o}^{\delta^o(G)}k_{\beta^{-1}},$$
where $\delta^o:G\to \sur{G}\times G$ sends $x$ to $(\sur{x},x)$, and $K^o=\{(\sur{a},a)\mid a\in K\}$. \par
Now we compute the tensor product $U\otimes_{k\sur{G}}V$ using Theorem 1.1 of~\cite{tensorbimodules}. Since $p_2\big(\delta(G)\big)=\sur{G}$, there is a single double coset $p_2\big(\delta(G)\big)\dom\sur{G}/p_1\big(\delta^o(G)\big)$. Moreover $k_2\big(\delta(G)\big)=1$, so we have
\begin{moneq}\label{U V}
U\otimes_{k\sur{G}}V\cong\Ind_{\delta(G)*\delta^o(G)}^{G\times G}\big(\Inf_{K}^{\delta(G)}k_\beta\otimes_k\Inf_{K^o}^{\delta^o(G)}k_{\beta^{-1}}\big),
\end{moneq}
where
\begin{align*}
\delta(G)*\delta^o(G)&=\big\{(a,b)\in G\times G\mid \exists \sur{c}\in\sur{G},\;(a,\sur{c})\in \delta(G)\,\hbox{and}\,(\sur{c},b)\in \delta^o(G)\big\}\\
&=\big\{(a,b)\in G\times G\mid \sur{a}=\sur{b}\big\}\\
&=\big\{(a,b)\in G\times G\mid ab^{-1}\in C_K(P)\big\}.
\end{align*}
Now if $a=xu$ and $b=yv$, with $x,y\in P$ and $u,v\in K$, we have
$$ab^{-1}=xuv^{-1}y^{-1}=x\cdot{^{uv^{-1}}(y^{-1})}\cdot uv^{-1},$$
so $ab^{-1}\in C_K(P)$ if and only if $uv^{-1}\in C_K(P)$ and $x=y$, i.e. $(u,v)\in\widetilde{K}$ and $(x,y)\in\Delta(P)$. It follows that $\delta(G)*\delta^o(G)=\Delta(P)\widetilde{K}$.\par
The action of $(a,b)\in \delta(G)*\delta^o(G)$ on the tensor product 
$$T:=\Inf_{K}^{\delta(G)}k_\beta\otimes_k\Inf_{K^o}^{\delta^o(G)}k_{\beta^{-1}}$$ is obtained as follows: Let $\sur{c}\in \sur{G}$ such that $(a,\sur{c})\in \delta(G)$ and $(\sur{c},b)\in \delta^o(G)$, that is $\sur{c}=\sur{a}=\sur{b}$. Then for $v\in \Inf_{K}^{\delta(G)}k_\beta$ and $w\in \Inf_{K^o}^{\delta^o(G)}k_{\beta^{-1}}$, we have $(a,b)\cdot (v\otimes w)=(a,\sur{c})\cdot v\otimes (\sur{c},b)\cdot w$. Here $T$ is one dimensional, with basis $1\otimes 1$, and
\begin{align*}
(a,b)\cdot(1\otimes 1)&=\beta\big(\varpi(a)\big)\beta\big(\varpi(b)\big)^{-1}(1\otimes 1)\\
&=\beta\big(\varpi(ab^{-1})\big)(1\otimes 1)\\
&=\lambda\big(\varpi(ab^{-1})\big)(1\otimes 1),
\end{align*}
since the restriction of $\beta$ to $C_K(P)$ is equal to $\lambda$. \smallskip\par
Now for $(a,b)\in \delta(G)*\delta^o(G)=\Delta(P)\widetilde{K}$, we have $\lambda\big(\varpi(ab^{-1})\big)=\widetilde{\lambda}\big(\rho(a,b)\big)$, where $\rho:\Delta(P)\widetilde{K}\to\widetilde{K}$ is the projection map. Then $T\cong\Inf_{\widetilde{K}}^{\Delta(P)\widetilde{K}}k_{\widetilde{\lambda}}$, and by~\ref{kGe} and~\ref{U V}, we get that 
$$U\otimes_{k\sur{G}}V\cong \Ind_{\delta(G)*\delta^o(G)}^{G\times G}T\cong \Ind_{\Delta(P)\widetilde{K}}^{G\times G}\Inf_{\widetilde{K}}^{\Delta(P)\widetilde{K}}k_{\widetilde{\lambda}}\cong kGe,$$
as was to be shown.\par
So if $C_K(P)\neq 1$, all the bimodules $kGe$ are mapped to 0 in $\CE_\SR(G)$, and the image of their direct sum $kG$ is also 0. Now as its identity element is equal to~0, the algebra $\CE_\SR(G)$ itself is equal to 0.\endpf
\begin{enonce}{Corollary} \label{K cyclic faithful}Let $G=P\rtimes K$, where $K$ is an elementary $p'$-group of order invertible in~$\SR$. If $\CE_\SR(G)\neq 0$, then $K$ is cyclic and acts faithfully on~$P$.
\end{enonce}
\pf Indeed if $\CE_\SR(G)\neq 0$, we know by Theorem~\ref{K cyclic} that $K$ is cyclic, and by Theorem~\ref{K faithful}, that $K$ acts faithfully on $P$.
\endpf
In other words $G=P\rtimes\langle u\rangle$, where $(P,u)$ is a $D^\Delta$-pair, as defined hereafter:\footnote{$D^\Delta$-pairs were first introduced in the slightly different Definition 4.4 of \cite{bouc-yilmaz}. The subsequent Lemma 4.5 there showed that the two definitions are equivalent.}
\vspace{-1ex}
\begin{enonce}{Definition} \label{D Delta pair}A $D^\Delta$-pair is a pair $(P,u)$ of a finite $p$-group $P$ and a $p'$-automorphism $u$ of $P$.
\end{enonce}
We recall the following notation (\cite{functorialequivalence}, Notation 6.8):
\begin{enonce}{Notation}  \label{Out(L,u)} For a $D^\Delta$-pair $(L,u)$, we denote by $\Aut(L,u)$ the group of automorphisms of the semidirect product $\Lu=L\rtimes \langle u\rangle$ which send $u$ to a conjugate of $u$, and by $\Out(L,u)$ the quotient $\Aut(L,u)/\Inn(\Lu)$ of this group by the group of inner automorphisms of $\Lu$.
\end{enonce}
\vspace{-2ex}
\section{Generators and relations for $\CE_\SR(G)$}
\begin{enonce}{Notation} For a finite group $G$, denote by $G^\natural$ the group $\Hom(G,k^\times)$. For $\lambda\in G^\natural$, let $k_\lambda$ denote the corresponding one dimensional $kG$-module. For $\gamma\in\Aut(G)$, denote\footnote{This generalizes Definition~4.1.3 of~\cite{ducellier}, up to replacing $\lambda$ with $\lambda^{-1}$, which is more convenient in our setting.} by $kG_{\gamma,\lambda}$ the $(kG,kG)$-bimodule equal to $kG$ as a vector space, with action given by
$$\forall (x,y,z)\in G^3,\;x\cdot z\cdot y=\lambda(y)^{-1}xz\gamma(y).$$
\end{enonce} 
\begin{enonce}{Lemma} \label{kG gamma lambda}Let $G$ be a finite group. 
\begin{enumerate}
\item Let $\gamma\in\Aut(G)$ and $\lambda\in G^\natural$. Then $kG_{\gamma,\lambda}$ is a diagonal $p$-permutation bimodule, and there is an isomorphism of $(kG,kG)$-bimodules
$$kG_{\gamma,\lambda}\cong \Ind_{\Delta_\gamma(G)}^{G\times G}k_{\lambda}.$$
\item \footnote{Up to the previous change of notation, this is Proposition 4.1.4 in~\cite{ducellier}.} Let $\delta\in\Aut(G)$ and $\mu\in G^\natural$. Then there is an isomorphism of $(kG,kG)$-bimodules
$$kG_{\gamma,\lambda}\otimes_{kG}kG_{\delta,\mu}\cong kG_{\gamma\circ\delta,(\lambda\circ\delta)\times \mu}.$$
\item If $kG$ has only one block, then $kG_{\gamma,\lambda}$ is an indecomposable $(kG,kG)$-bimodule, for any $\gamma\in\Aut(G)$ and any $\lambda\in G^\natural$.\vspace{-2ex}
\end{enumerate}
\end{enonce}
\pf 1. Let $S$ be a Sylow $p$-subgroup of $G$. Then $\lambda(x)=1$ for any $x\in S$, so the restriction of $kG_{\gamma,\lambda}$ to the Sylow $p$-subgroup $S\times S$ of $G\times G$ is the permutation bimodule $kG$, with action $x\cdot z\cdot y=xz\gamma(y)$, for $x,y\in S$ and $z\in G$. Moreover this action is free on both sides, so $kG_{\gamma,\lambda}$ is a diagonal $p$-permutation bimodule. Finally, the map $g\in G\mapsto (g,1)\Delta_\gamma(G)$ is a bijection from $G$ to $(G\times G)/\Delta_\gamma(G)$, and using this bijection, it is easy to check that $kG_{\gamma,\lambda}\cong \Ind_{\Delta_\gamma(G)}^{G\times G}k_{\lambda}$.\mpn
2. The map $(g\otimes g')\mapsto \lambda(g')^{-1}g\gamma(g')$, for $g,g'\in G$, from $kG_{\gamma,\lambda}\otimes_{kG}kG_{\delta,\mu}$ to $kG_{\gamma\circ\delta,(\lambda\circ\delta)\times \mu}$ induces a well defined isomorphism of $(kG,kG)$-bimodules.\mpn
3. It is clear that for any $(kG,kG)$-bimodule $M$, any $\delta\in\Aut(G)$, and any $\mu\in G^\natural$, the $k$-vector space $M\otimes_{kG} kG_{\delta,\mu}$ is isomorphic to $M$. So if $kG_{\gamma,\lambda}$ splits as a direct sum of non-zero $(kG,kG)$-bimodules $M$ and $M'$, the tensor product $kG_{\gamma,\lambda}\otimes_{kG} kG_{\delta,\mu}$, splits as the direct sum of $(kG,kG)$-bimodules $M\otimes kG_{\delta,\mu}$ and $M'\otimes kG_{\delta,\mu}$, none of which is equal to zero. Then $kG_{\gamma\circ\delta,(\lambda\circ \delta)\times \mu}$ splits as a direct sum of non-zero bimodules. Taking $\delta=\gamma^{-1}$ and $\mu=(\lambda\circ\gamma^{-1})^{-1}$, we get that the $(kG,kG)$-bimodule $kG_{\Id,1}\cong kG$ splits non-trivially. So $kG$ is not indecomposable, that is, $kG$ has more than one block.
\endpf
In the rest of this section, in view of Corollary~\ref{K cyclic faithful}, we assume the following:
\begin{enonce}{Hypothesis} \label{hyp}The group $G$ is of the form $P\rtimes K$, where $P$ is a $p$-group and $K$ is a cyclic $p'$-group of order invertible in $\SR$, acting faithfully on $P$.
\end{enonce}
{\flushleft We} want to find the structure of the algebra $\CE_\SR(G)$. 
First we look for generators of $\CE_\SR(G)$ as an $\SR$-module. 

\begin{enonce}{Lemma} \label{essential bimodule}Assume that~\ref{hyp} holds. Then:
\begin{enumerate}
\item $kG$ has only one block. So for $\gamma\in\Aut(G)$ and $\lambda\in G^\natural$, the bimodule $kG_{\gamma,\lambda}$ is indecomposable, with vertex $\Delta_\gamma(P)$. 
\item Let $\gamma, \delta\in\Aut(G)$ and $\lambda,\mu\in G^\natural$. Then the bimodules $kG_{\gamma,\lambda}$ and $kG_{\delta,\mu}$ are isomorphic if and only if $\lambda=\mu$ and $\delta\circ\gamma^{-1}$ is an inner automorphism of $G$.
\item Conversely, if $M$ is an indecomposable diagonal $p$-permutation bimodule with vertex $\Delta_\gamma(P)$, for $\gamma\in\Aut(G)$, then there exists $\lambda\in G^\natural$ such that $M\cong kG_{\gamma,\lambda}$.
\item In particular, if $M$ is an essential indecomposable diagonal $p$-permuta\-tion $(kG,kG)$-bimodule, there exist $\gamma\in \Aut(G)$ and $\lambda\in G^\natural$ such that $M\cong kG_{\gamma,\lambda}$ as $(kG,kG)$-bimodule. \end{enumerate}
\end{enonce}
\pf 1. The $p$-subgroup $P$ of $G$ is normal, and $C_G(P)=Z(P)\times C_K(P)=Z(P)\leq P$ since $C_K(P)=1$ as $K$ acts faithfully on $P$. So $kG$ has only one block (\cite{benson1} Proposition 6.2.2). Then all the bimodules $kG_{\gamma,\lambda}$ are indecomposable, by Lemma~\ref{kG gamma lambda}.\par
Let $\gamma\in\Aut(G)$. Set $N:=N_{G\times G}\big(\Delta_\gamma(P)\big)$, and $\sur{N}=N/\Delta_\gamma(P)$. Then 
$$N=\{(a,b)\in G\times G\mid a^{-1}\gamma(b)\in C_G(P)=Z(P)\},$$
so the second projection $p_2$ induces a short exact sequence
$$\xymatrix{
1\ar[r]&Z(P)\ar[r]&\sur{N}\ar[r]&K\ar[r]&1,
}
$$
which is split (by the map $x\in K\mapsto \big(\gamma(x),x\big)\Delta_\gamma(P)$). Now the indecomposable projective $k\sur{N}$-modules are the modules $\Ind_{K}^{\sur{N}}k_\lambda$, for $\lambda\in K^\natural$. Moreover
\begin{align*}
\Ind_{N}^{G\times G}\Inf_{\sur{N}}^N\Ind_{K}^{\sur{N}}k_\lambda&\cong \Ind_{N}^{G\times G}\Ind_{\Delta_\gamma(P)\cdot\Delta_\gamma(K)}^Nk_\lambda\\
&\cong\Ind_{\Delta_\gamma(G)}^{G\times G}k_\lambda\cong kG_{\gamma,\lambda},
\end{align*}
by Lemma~\ref{kG gamma lambda}.  This gives another proof that $kG_{\gamma,\lambda}$ is indecomposable, and also shows that it has vertex $\Delta_\gamma(P)$, and Brauer quotient
\begin{moneq}\label{brauer kG gamma lambda}
kG_{\gamma,\lambda}[\Delta_\gamma(P)]\cong \Ind_K^{\sur{N}}k_\lambda\cong kZ(P)\otimes_kk_{\lambda}.
\end{moneq}
{\flushleft{2.}} First, if $\lambda=\mu$ and $\delta=i_x\circ \gamma$, where $i_x$ is conjugation by $x\in G$, then one checks easily that the map $g\in G\mapsto gx\in G$ induces a bimodule isomorphism $kG_{\gamma,\lambda}\cong kG_{\delta,\mu}$. \par
For the converse, let $\gamma'=\gamma^{-1}$ and $\lambda'=(\lambda\circ\gamma')^{-1}$. Then by Lemma~\ref{kG gamma lambda}, we have an isomorphism of $(kG,kG)$-bimodules
$$kG_{\gamma,\lambda}\otimes_{kG}kG_{\gamma',\lambda'}\cong kG_{\gamma\circ\gamma',(\lambda\circ \gamma')\times \lambda'}=kG_{\Id,1}=kG.$$
So if $kG_{\delta,\mu}\cong kG_{\gamma,\lambda}$, we have an isomorphism of $(kG,kG)$-bimodules
$$kG_{\delta,\mu}\otimes_{kG}kG_{\gamma',\lambda'}\cong kG,$$
that is
$$kG_{\delta\circ\gamma',(\mu\circ\gamma')\times \lambda'}\cong kG_{\Id,1}.$$
So if we know that an isomorphism of bimodules $kG_{\theta,\rho}\cong kG_{\Id,1}$, where $\theta\in\Aut(G)$ and $\rho\in G^\natural$ implies that $\theta$ is inner and $\rho=1$, we are done, since we can conclude that $\delta\circ\gamma'=\delta\circ \gamma^{-1}$ is inner, and that $(\mu\circ\gamma')\times \lambda'=1$, i.e. $\mu=\lambda$. In other words, we can assume $\gamma=\Id$ and $\lambda=1$, and  that $kG_{\delta,\mu}\cong kG$. \par
 Now if $kG_{\delta,\mu}$ is isomorphic to $kG$, then its vertex $\Delta_\delta(P)$ is contained in -~hence equal to~- $\Delta(P)$, up to conjugation in $G\times G$. It means that the restriction of $\delta$ to $P$ is equal to the conjugation by some element of $G$. Up to composing $\delta$ with some inner automorphism of $G$, we can assume that this restriction is equal to the identity, and then $\delta$ is equal to the conjugation by some element of $Z(P)$, by Lemma~\ref{inner}. This shows that $\delta$ is inner. \par
Then we have a bimodule isomorphism $kG_{\delta,\lambda}\cong kG_{\Id,\lambda}$ by the first remark in the proof of Assertion 2. In other words, we can assume $\delta=\Id$ and $kG_{\Id,\lambda}\cong kG_{\Id,1}$. Then the Brauer quotients at $\Delta(P)$ of these bimodules are isomorphic. Hence $kZ(P)\otimes k_\lambda\cong kZ(P)\otimes k\cong kZ(P)$. Now the fixed points of $Z(P)$ on $kZ(P)\otimes k_\lambda$ form a $kK$-module isomorphic to $k_\lambda$, so $k_\lambda\cong k$ as $K$-module, hence $\lambda=1$, as was to be shown.
 \mpn
3. Suppose conversely that $M$ is an indecomposable diagonal $p$-permutation $(kG,kG)$-bimodule with vertex $\Delta_\gamma(P)$, where $\gamma\in\Aut(G)$. Then $M[\Delta_\gamma(P)]$ is an indecomposable projective $k\sur{N}$-module, of the form $\Ind_K^{\sur{N}}k_\lambda$ for some $\lambda\in K^\natural$, and then
$$M\cong\Ind_N^{G\times G}\Inf_{\sur{N}}^N\Ind_K^{\sur{N}}k_\lambda\cong kG_{\gamma,\lambda}.$$
4. As in the proof of Theorem~\ref{K cyclic}, we apply Lemma~\ref{does not factor}, in the case $H=G$ and $U=M$. We know that
$$M\cong \Ind_{\Delta_\pi(P)\cdot T}^{G\times G}\Inf_{T}^{\Delta_\varphi(P)\cdot T}W,$$
for some $\pi\in \Aut(P)$, some subgroup $T$ of $N_{K\times K}\big(\Delta_\pi(P)\big)$ with $p_2(T)=K$, and some simple $kT$-module $W$. \par
Moreover $k_1(T)\leq k_1\Big(N_{G\times G}\big(\Delta_\pi(P)\big)\Big)=C_G(P)=Z(P)$. So $k_1(T)=1$ since $T$ is a $p'$-group. Then $p_2:T\to K$ is an isomorphism, with inverse $\theta$, and $T=\Delta_\theta(K)$. Then $T$ is cyclic, and $W\cong k_\lambda$ for some $\lambda\in T^\natural$.\smallskip\par
Now $T$ normalizes $\Delta_\pi(P)$ if and only if $^{\theta(x)}\pi(y)=\pi(^xy)$ for any $x\in K$ and any $y\in P$. Then the map $\gamma:y\cdot x\mapsto \pi(y)\cdot\theta(x)$, where $y\in P$ and $x\in K$, is an automorphism of $G$, such that $\gamma(K)=K$, and $\Delta_\pi(P)\cdot L=\Delta_\gamma(G)$. Then 
$$M\cong\Ind_{\Delta_\gamma(G)}^{G\times G}k_\mu,$$
where $\mu=\lambda\circ \varpi\in G^\natural\cong \Delta_\gamma(G)^\natural$. Now $M\cong kG_{\gamma,\mu}$ by Lemma~\ref{kG gamma lambda}.\medskip \endpf

 It follows that $\CE_\SR(G)$ is linearly generated by the images of the $(kG,kG)$-bimodules $kG_{\gamma,\lambda}$, for $\gamma\in\Aut(G)$ and $\lambda\in G^\natural$. By Lemma~\ref{essential bimodule}, we can take $\gamma$ in a set of representatives of elements of $\Out(G)$. We want to describe the linear relations between these generators. In other words, we want to find equalities of the form
\begin{moneq}\label{linear combination}
\sum_{\substack{\gamma\in\Out(G)\\\lambda\in G^\natural}}r_{\gamma,\lambda}\,kG_{\gamma,\lambda}=\sum_{i=1}^n s_i\, U_i\otimes_{kH_i}V_i,
\end{moneq} 
{\flushleft in} $\SR T^\Delta(G,G)$, where $r_{\gamma,\lambda}\in\SR$, $n\in\NN$, and for $1\leq i\leq n$, $H_i$ is a finite group of order smaller than the order of $G$, $U_i$ is a diagonal $p$-permutation $(kG,kH_i)$-bimodule, $V_i$ is a diagonal $p$-permutation $(kH_i,kG)$-bimodule, and $s_i\in \SR$. We can assume moreover that for $1\leq i\leq n$, the essential algebra $\CE_\SR(H_i)$ is non-zero: Indeed otherwise, the identity bimodule $kH_i\in\SR T^\Delta(H_i,H_i)$ is a linear combination with coefficients in $\SR$ of elements of $\SR T^\Delta(H_i,X)\otimes_{kX}\SR T^\Delta(X,H_i)$, for $|X|<|H_i|$, and we can replace $H_i$ by smaller groups in~(\ref{linear combination}). \par
 
Hence, by Corollary~\ref{semidirect}, we can assume that for $1\leq i\leq n$, we have $H_i=Q_i\rtimes L_i$, where $Q_i$ is a $p$-group, and $L_i$ is an elementary $p'$-group. We can also assume that $U_i$ and $V_i$ are indecomposable, and that $U_i$ is right essential and $V_i$ is left essential. 
\begin{enonce}{Lemma} \label{H in G}Assume that~\ref{hyp} holds.  Let $H$ be a finite group, and $U$ be a right essential indecomposable diagonal $p$-permutation $(kG,kH)$-bimodule. Then the vertices of $U$ have order at most $|P|$. If $U$ has vertex of order $|P|$, then there exists an injective group homomorphism $\sigma: H\hookrightarrow G$ such that $P\leq \sigma(H)$ and $\lambda\in H^\natural=\Delta_\sigma(H)^\natural$ such that
$$U\cong\Ind_{\Delta_\sigma(H)}^{G\times H}k_\lambda.$$
\end{enonce}
\pf We know from Lemma~\ref{does not factor} that $H=Q\rtimes L$, where $Q$ is a $p$-group with an embedding $\pi:Q\hookrightarrow P$, and $L$ is an elementary $p'$-group. Moreover there is a subgroup $T$ of $N_{K\times L}\big(\Delta_\pi(Q)\big)$ with $p_2(T)=L$, and a simple $kT$-module $W$ such that
$$U\cong\Ind_{\Delta_\pi(Q)\cdot T}^{G\times H}\Inf_{T}^{\Delta_\pi(Q)\cdot T}W.$$
Then a vertex of $U$ is contained in $\Delta_\pi(Q)\cdot T$ up to conjugation, so it has order at most $|Q|\leq |P|$. And if it has order $|P|$, the embedding $\pi:Q\hookrightarrow P$ is an isomorphism. Moreover $k_1(T)\leq k_1\Big(N_{G\times H}\big(\Delta_\pi(P)\big)\Big)=C_G(P)$, so $k_1(T)\leq C_K(P)=1$ since $K$ acts faithfully on~$P$. Then the projection map $p_2:T\to L$ is an isomorphism, and then $T=\Delta_\tau(L)$, for some injective group homomorphism $\tau:L\hookrightarrow K$. In particular $L$ and $T$ are cyclic. \par
Moreover since $T=\Delta_\tau(L)\leq N_{G\times H}\big(\Delta_\pi(P)\big)$, we have $^{\tau(l)}\pi(x)=\pi({^lx})$ for any $l\in L$ and $x\in Q$. Then the map $\sigma:H=Q\cdot L\to G$ sending $x\cdot l$, for $x\in Q$ and $l\in L$, to $\pi(x)\cdot\tau(l)$, is an injective group homomorphism, such that $P\leq \gamma(H)\leq G$. Moreover $\Delta_\pi(Q)\cdot T=\Delta_\sigma(H)$. Finally $T$ is cyclic, so there is $\lambda\in T^\natural=L^\natural\cong H^\natural$ such that $W=k_{\lambda}$, and $U\cong\Ind_{\Delta_\sigma(H)}^{G\times H}k_\lambda$, completing the proof.\endpf
\begin{enonce}{Theorem} \label{induced}Assume that~\ref{hyp} holds. Let $H$ be a finite group, let $U$ (resp. $V$) be a right (resp. left) essential indecomposable diagonal $p$-permutation $(kG,kH)$-bimodule (resp. $(kH,kG)$-bimodule).
\begin{enumerate}
\item If $U\otimes_{kH}V$ has an indecomposable direct summand with vertex of order $|P|$, then there is a subgroup $I\cong H$ of $G$, containing $P$, an automorphism $\psi$ of $I$, and $\zeta\in I^\natural=\Delta_\psi(I)^\natural$ such that
$$U\otimes_{kH}V\cong \Ind_{\Delta_\psi(I)}^{G\times G}k_\zeta.$$
\item If $U\otimes_{kH}V$ admits an essential indecomposable summand, then this summand has vertex $\Delta_\gamma(P)$ for some $\gamma\in\Aut(G)$, and there exists $J\leq K$ and $\theta\in J^\natural$ such that $P\cdot J\cong H$ and
$$U\otimes_{kH}V\cong\bigoplus_{\alpha\in \CI_\theta}kG_{\gamma,\alpha}$$
as $(kG,kG)$-bimodules, where $\CI_\theta=\big\{\alpha\in G^\natural=K^\natural\mid \Res_J^K\alpha=\theta\big\}$. 
\end{enumerate}
\end{enonce}
\pf 1. Let $X$ be a vertex of $U$. Then $X$ is a diagonal subgroup of $P\times H$, so $|X|\leq |P|$, and $U$ is a direct summand of $\Ind_{X}^{G\times H}k$. Similarly, if $Y$ is a vertex of $V$, then $Y$ is a diagonal subgroup of $H\times P$, hence $|Y|\leq |P|$, and $V$ is a direct summand of $\Ind_{Y}^{H\times G}k$. It follows that $U\otimes_{kH}V$ is a direct summand of
$$\bigoplus_{h\in p_2(X)\dom H/p_1(Y)}\Ind_{X*{^{(h,1)}Y}}^{G\times G}k.$$
So the vertices of the indecomposable summands of $U\otimes_{kH}V$ are contained (up to conjugation) in some group $X*{^{(h,1)}Y}$, which has order at most $\min(|X|,|Y|)$. If one of them has order $|P|$, then $|P|\leq |X|\leq|P|$, hence $|X|=|P|$, and $|P|\leq |Y|\leq |P|$, so $|Y|=|P|$. \par
By Lemma~\ref{H in G}, it follows that there is an embedding $\sigma:H\hookrightarrow G$ with $P\leq \sigma(H)\leq G$, and $\lambda\in H^\natural$, such that $U\cong \Ind_{\Delta_\sigma(H)}^{G\times H}k_\lambda$. Similarly, swapping $G$ and $H$, there is an embedding $\tau:H\to G$ with $P\leq\tau(H)\leq G$, and $\mu\in H^\natural$, such that $V\cong \Ind_{\Delta_\tau^o(H)}^{H\times G}k_\mu$, where $\Delta_\tau^o(H)=\Big\{\big(h,\tau(h)\big)\mid h\in H\Big\}$. Then
$$U\otimes_{kH}V\cong \Ind_{\Delta_\sigma(H)*\Delta_\tau^o(H)}^{G\times G}(k_\lambda\otimes_kk_\mu).$$
Now $\sigma(H)/P$ and $\tau(H)/P$ are subgroups of the same order of the cyclic group $K$, so $\sigma(H)=\tau(H)$. Set $I:=\tau(H)$. There is a unique automorphism $\psi$ of $I$ such that $\psi\big(\tau(h)\big)=\sigma(h)$ for all $h\in H$. Then
$$\Delta_\sigma(H)*\Delta_\tau^o(H)=\Big\{\big(\sigma(h),\tau(h)\big)\mid h\in H\Big\}=\big\{(\psi(x),x)\mid x\in I\big\}=\Delta_\psi(I).$$
Moreover $k_\lambda\otimes k_\mu$ is one dimensional, so there is a unique $\zeta\in I^\natural$ such that $k_\lambda\otimes k_\mu\cong k_\zeta$ as $kI$-modules, defined by $\zeta(x)=(\lambda\mu)\big(\tau^{-1}(x)\big)$ for $x\in I$. This completes the proof of Assertion 1.\mpn
2. By Lemma~\ref{essential bimodule}, an essential diagonal $p$-permutation bimodule $M$ is isomorphic to $kG_{\gamma,\lambda}$ for some $\gamma\in\Aut(G)$ and $\lambda\in G^\natural$. Then $M$ has vertex $\Delta_\gamma(P)$, of order $|P|$, so the conclusion of Assertion 1 holds. Hence there is a subgroup $I$ of $G$, containing $P$, an automorphism $\psi$ of $I$, and $\zeta\in I^\natural$, such that

$$U\otimes_{kH} V\cong\Ind_{\Delta_\psi(I)}^{G\times G}k_\zeta.$$
In particular $\Delta_\gamma(P)$ is contained in $\Delta_\psi(P)$, up to conjugation, and we can assume that $\Delta_\gamma(P)=\Delta_\psi(P)$, i.e. that $\psi$ is the restriction of $\gamma$ to $P$. We have $\Delta_\psi(P)\leq \Delta_\psi(I)\leq N:=N_{G\times G}\big(\Delta_\psi(P)\big)=N_{G\times G}\big(\Delta_\gamma(P)\big)$. So $N$ fits in a short exact sequence of groups
$$\xymatrix{
1\ar[r]&\ar[r]Z(P)\ar[r]& N\ar[r]^-{p_2}&G\ar[r]&1.
}
$$
Similarly, the group $\sur{N}:=N/\Delta_\psi(P)$ fits in the sequence
$$\xymatrix{
1\ar[r]&\ar[r]Z(P)\ar[r]& \sur{N}\ar[r]&K\ar[r]&1.
}
$$
This sequence splits via $x\in K\mapsto \big(\gamma(x),x\big)\Delta_\psi(P)$, and we view $K$ as a subgroup of $\sur{N}$ via this map.\par
The group $\sur{\Delta}_\psi(I)=\Delta_\psi(I)/\Delta_\psi(P)$ is a subgroup of $\sur{N}$, and intersects $Z(P)$ trivially. So $\sur{\Delta}_\psi(I)$ is isomorphic to a subgroup $J$ of $K$. Let $\theta:J\to k^\times$ be the image of $\sur{\zeta}$ under this isomorphism.\par
Since $\Delta_\psi(P)$ acts trivially on $k_\zeta$, we have $k_\zeta=\Inf_{\sur{\Delta}_\psi(I)}^{\Delta_\psi(I)}k_{\sur{\zeta}}$, where $\sur{\zeta}$ is the homomorphism $\sur{\Delta}_\psi(I)\to k^\times$ corresponding to $\zeta$. Hence we have
\begin{align*}
U\otimes_{kH}V&\cong \Ind_{\Delta_\psi(I)}^{G\times G}k_\zeta\\
&=\Ind_{N}^{G\times G}\Ind_{\Delta_\psi(I)}^N\Inf_{\sur{\Delta}_\psi(I)}^{\Delta_\psi(I)}k_{\sur{\zeta}}\\
&\cong \Ind_{N}^{G\times G}\Inf_{\sur{N}}^N\,\Ind_{\sur{\Delta}_\psi(I)}^{\sur{N}}k_{\sur{\zeta}}\\
&\cong\Ind_{N}^{G\times G}\Inf_{\sur{N}}^N\,\Ind_{\sur{\Delta}_\psi(I)}^{\sur{N}}\Iso_{J}^{\sur{\Delta}_\psi(I)}k_\theta\\
&\cong \Ind_{N}^{G\times G}\Inf_{\sur{N}}^N\Ind_K^{\sur{N}}\Ind_J^Kk_\theta.
\end{align*}
Now $K$ is cyclic, so $\Ind_J^Kk_\theta=\mathop{\sum}_{\substack{\alpha\in K^\natural\\\Res_J^K\alpha=\theta}}\limits k_\alpha$, and if $\alpha\in K^\natural$, then $\Ind_K^{\sur{N}}k_\alpha$ is an indecomposable projective $k\sur{N}$-module. Then $\Ind_{N}^{G\times G}\Inf_{\sur{N}}^N\Ind_K^{\sur{N}}k_\alpha\cong kG_{\gamma,\alpha}$. 
Now
$$U\otimes_{kH}V\cong\bigoplus_{\beta\in \CI_\theta}kG_{\gamma,\beta},$$
as was to be shown.
\endpf
\begin{enonce}{Notation} Assume that~\ref{hyp} holds.
\begin{enumerate}
\item We abuse notation identifying $\lambda\in K^\natural$ with $k_\lambda\in R_k(K)$.
\item We set $\sur{R}_k(K)=R_k(K)/\sum_{L<K}\limits\Ind_L^KR_k(L)$, and we let $\alpha\mapsto \sur{\alpha}$ denote the projection map.
\item Let $\gamma\in \Aut(G)$. Then $\sur{N}_{G\times G}\big(\Delta_\gamma(P)\big)\cong Z(P)\rtimes K$, so taking coinvariants by $Z(P)$ yields an isomorphism 
$$v\in \mathrm{Proj}\Big(k\sur{N}_{G\times G}\big(\Delta_\gamma(P)\big)\Big)\mapsto v_{Z(P)}\in R_k(K).$$
For $u\in T^\Delta(G,G)$, let $r_\gamma(u)$ denote $\sur{u[\Delta_\gamma(P)]_{Z(P)}}\in \sur{R}_k(K)$.
\end{enumerate}
\end{enonce}
{\flushleft We} note that $\sum_{L<K}\limits\Ind_L^KR_k(L)$ is an ideal of the ring $R_k(K)$, so $\sur{R}_k(K)$ has a natural quotient ring structure. Moreover, the group $\Aut(G)$ acts on $G^\natural=K^\natural$, and $\Inn(G)$ acts trivially on $G^\natural$, since $[G,G]\leq P$. So $\Out(G)$ acts on $R_k(K)$ and $\sur{R}_k(K)$ by ring automorphisms. 
\begin{enonce}{Notation} \label{ring semidirect}We denote by $\Out(G)\ltimes \sur{R}_k(K)$ the semidirect product of $\Out(G)$ with $\sur{R}_k(K)$, i.e.
$$\Out(G)\ltimes \sur{R}_k(K)=\bigoplus_{\gamma\in\Out(G)}\gamma\ltimes \sur{R}_k(K),$$
where $\gamma\ltimes \sur{R}_k(K)$ denotes a copy of $\sur{R}_k(K)$ indexed by $\gamma$. \par
Then $\Out(G)\ltimes \sur{R}_k(K)$ is a ring for the product defined by
$$\forall \gamma,\delta\in\Out(G),\,\forall \lambda,\mu\in R_k(K),\;(\gamma\ltimes \sur{\lambda})\cdot(\delta\ltimes\sur{\mu}):=(\gamma\circ\delta)\ltimes\sur{\big((\lambda\circ\delta)\times\mu\big)}.$$
We set 
$$\Out(G)\ltimes\SR \sur{R}_k(K):=\SR\otimes_\ZZ \big(\Out(G)\ltimes \sur{R}_k(K)\big)\cong \bigoplus_{\gamma\in\Out(G)}\gamma\ltimes \SR \sur{R}_k(K).$$
\end{enonce}

In the next statement, we recall explicitly Hypothesis~\ref{hyp}, for the reader's convenience.
\begin{enonce}{Theorem} \label{structure essential}Let $G$ be a group of the form $P\rtimes K$, where $P$ is a $p$-group and $K$ is a cyclic $p'$-group of order invertible in $\SR$, acting faithfully on $P$. Then:
\begin{enumerate}
\item The map
$$\gamma\ltimes \sur{\alpha}\in \Out(G)\ltimes \SR \sur{R}_k(K)\mapsto \Sur{kG_{\gamma,\alpha}}\in\CE_\SR(G),$$
where $\alpha\in K^\natural$ and $\gamma\in \Out(G)$, extends to a well defined algebra homomorphism $\St$.
\item The map
$$\Sur{u}\in \CE_\SR(G)\mapsto\sum_{\gamma\in\Out(G)}\gamma\ltimes r_\gamma(u)\in \Out(G)\ltimes \SR \sur{R}_k(K),$$
where $u\in \SR T^\Delta(G,G)$, is a well defined algebra homomorphism $\Ss$.
\item The maps
$$\xymatrix{\CE_\SR(G)\ar[r]<.5ex>^-{\Ss}&\Out(G)\ltimes \SR \sur{R}_k(K)\ar[l]<.5ex>^-{\St}
}
$$ are isomorphisms of algebras, inverse to each other.
\end{enumerate}
\end{enonce}
\pf Proving that the map $\St$ is well defined amounts to proving that if $u\in K^\natural$ is induced from a proper subgroup $J$ of $K$, and if $\gamma\in \Out(G)$, then $\St(\gamma\ltimes\sur{u})=0$. Let $J<K$, and $\theta\in J^\natural$. 
Then $u=\Ind_J^K\theta=\sum_{\alpha\in\CI_\theta}\limits k_\alpha$, so
$$\St(\gamma\ltimes u)=\CE_\SR\Big(\bigoplus_{\alpha\in \CI_\theta}kG_{\gamma,\alpha}\Big).$$
But setting $N:=N_{G\times G}\big(\Delta_\gamma(P)\big)$ and $\sur{N}:=N/\Delta_\gamma(P)$, we have
\begin{align*}
\bigoplus_{\alpha\in \CI_\theta}kG_{\gamma,\alpha}&\cong\bigoplus_{\alpha\in \CI_\theta}\Ind_{N}^{G\times G}\Inf_{\sur{N}}^N\Ind_K^{\sur{N}}k_\alpha\\
&\cong\Ind_{N}^{G\times G}\Inf_{\sur{N}}^N\Ind_K^{\sur{N}}\Big(\bigoplus_{\alpha\in \CI_\theta}k_\alpha\Big)\\
&\cong \Ind_{N}^{G\times G}\Inf_{\sur{N}}^N\Ind_K^{\sur{N}}\Ind_J^Kk_\theta\\
&\cong \Ind_{N}^{G\times G}\Ind_{\Delta_\gamma(P)\cdot\Delta_\gamma(K)}^N\Inf_K^
{\Delta_\gamma(P)\cdot\Delta_\gamma(K)}\Ind_J^Kk_\theta\\
&\cong \Ind_{\Delta_\gamma(G)}^{G\times G}\Inf_K^{\Delta_\gamma(G)}\Ind_J^Kk_\theta\\
&\cong \Ind_{\Delta_\gamma(G)}^{G\times G}\Ind_{\Delta_\gamma(P)\cdot\Delta_\gamma(J)}^{\Delta_\gamma(G)}\Inf_J^{\Delta_\gamma(P)\cdot\Delta_\gamma(J)}k_\theta\\
&\cong \Ind_{\Delta_\gamma(P\cdot J)}^{G\times G}\Inf_J^{\Delta_\gamma(P\cdot J)}k_\theta.
\end{align*}
But $p_2\big(\Delta_\gamma(P\cdot J)\big)=P\cdot J<G$ since $J<K$. Hence $\CE_\SR\Big(\bigoplus_{\alpha\in \CI_\theta}kG_{\gamma,\alpha}\Big)=0$, as was to be shown, so the map $\St$ is well defined.\par
Now comparing the products in Lemma~\ref{kG gamma lambda} and Notation~\ref{ring semidirect}, we get that $\St$ is a homomorphism of $\SR$-algebras. Moreover the identity element of $\Out(G)\ltimes \SR \sur{R}_k(K)$, which is $\Id\ltimes \sur{1}$, is mapped by $\St$ to $\CE_\SR(kG_{\Id,1})=\CE_\SR(kG)$, which is the identity element of $\CE_\SR(G)$.
\mpn
2. Proving that the map $\Ss$ is well defined amounts to proving that if $H$ is a finite group with $|H|<|G|$, if $U$ (resp. $V$) is a right (resp. left) essential diagonal $p$-permutation $(kG,kH)$-bimodule (resp. $(kH,kG)$-bimodule), then $\Ss(U\otimes_{kH}V)=0$. So let $\gamma\in \Aut(G)$ such that $r_\gamma(U\otimes_{kH}V)\neq 0$. Then in particular $(U\otimes_{kH}V)[\Delta_\gamma(P)]\neq 0$, so $U\otimes_{kH}V$ admits an indecomposable summand with vertex $\Delta_\gamma(P)$. By Lemma~\ref{essential bimodule}, this summand is isomorphic to $kG_{\gamma,\lambda}$, for some $\lambda\in G^\natural$. By Theorem~\ref{induced}, there is a subgroup $J$ of $K$ with $P\cdot J\cong H$, and $\theta\in J^\natural$ such that
$$U\otimes_{kH}V\cong\bigoplus_{\alpha\in \CI_\theta}kG_{\gamma,\alpha}.$$
Now by~(\ref{brauer kG gamma lambda})
$$(U\otimes_{kH}V)[\Delta_\gamma(P)]\cong \bigoplus_{\alpha\in\CI_\theta}\Ind_K^{\sur{N}}k_\alpha\cong\Ind_K^{\sur{N}}\big(\bigoplus_{\alpha\in I_\theta}k_\alpha\big)\cong kZ(P)\otimes_k\Ind_J^Kk_\theta.$$
It follows that in $\sur{R}_k(K)$, we have $r_\gamma(U\otimes_{kH}V)=\sur{\Ind_J^Kk_\theta}=0$ since $J<K$ as $P\cdot J\cong H$ and $|H|<|G|=|P||K|$. This contradiction shows that $\mathsf{S}$ is well defined.\par 
We postpone the proof that $\Ss$ is an algebra homomorphism at the end of the proof of Assertion 3.\mpn
3. Let $\gamma\in\Aut(G)$ and $\alpha\in K^\natural$. Then $kG_{\gamma,\alpha}$ has vertex $\Delta_\gamma(P)$, and $kG_{\gamma,\alpha}[\Delta_{\gamma}(P)]\cong kZ(P)\otimes k_\alpha$ by~(\ref{brauer kG gamma lambda}). We also get from Lemma~\ref{essential bimodule} that $kG_{\gamma,\alpha}[\Delta_{\delta}(P)]=0$ if $\delta\in \Aut(G)$ and $\delta\neq \gamma$ in $\Out(G)$: Indeed, if $\Delta_\gamma(P)$ and $\Delta_\delta(P)$ are conjugate in $G\times G$, then the restriction of $\delta$ to $P$ is equal to the restriction of $\delta$ to $P$, up to an automorphism given by conjugation by some element of $G$, which we may assume to be trivial. Then $\delta^{-1}\gamma$ is inner, by Lemma~\ref{inner}. \par
It follows that $\Ss\big(\CE_\SR(kG_{\gamma,\alpha)})\big)=\gamma\ltimes \sur{\alpha}$. Since $\St(\gamma\ltimes \sur{\alpha})=\CE_\SR(kG_{\gamma,\alpha})$, the maps $\Ss$ and $\St$ are inverse to each other. In particular, they are bijections, so $\Ss$ is a map of $\SR$-algebras, as $\St$ is. This completes the proof.\endpf
\section{The simple functors}
\tpar We want to consider the simple diagonal $p$-permutation functors, so by general arguments, we can assume that our ring $\SR$ of coefficients is a field~$\F$. Moreover, in order to apply the results of the previous sections, we want that $p'$-groups have order invertible in $\F$. So we are left with the cases where $\F$ has characteristic 0 or $p$.\par
If $S$ is a simple diagonal $p$-permutation functor over $\F$, and $H$ is a finite group of minimal order such that $S(H)\neq 0$, then $V=S(H)$ is a simple module for the essential algebra $\CE_\F(H)$. In particular $\CE_\F(H)\neq 0$, so $H=\Lu$ for some $D^\Delta$-pair $(L,u)$. Moreover, we have an isomorphism of algebras
$$\mathcal{E}_\FF\big(\Lu\big)\cong \Out\big(\Lu\big)\ltimes \F\sur{R}_k(\Lu).\vspace{-2ex}$$
\tpar Conversely, if $(L,u)$ is a $D^\Delta$-pair, and $V$ is a simple $\CE_\F\big(\Lu\big)$-module, then we denote by $S_{\Lu,V}$ the unique simple diagonal $p$-permutation functor with minimal group $\Lu$ and such that $S_{\Lu,V}\big(\Lu\big)\cong V$ as an $\CE_\F\big(\Lu\big)$-module. 
The evaluation of $S_{\Lu,V}$ at a finite group $G$ is isomorphic to 
$$S_{\Lu,V}(G)\cong \Big(\FF T^\Delta\big(G,\Lu\big)\otimes_{\mathcal{E}_\FF(\Lu)}V\Big)/\mathcal{R},$$
where $\mathcal{R}$ is the subspace generated by all finite sums $\sum_{i\in I}f_i\otimes v_i$, with $f_i\in \FF T^\Delta\big(G,\Lu\big)$ and $v_i\in V$ for $i\in I$, such that $\sum_{i\in I}\pi(\varphi\circ f_i)\cdot v_i=0$ for all $\varphi\in \FF T^\Delta(\Lu,G)$, where $\pi: \FF T^\Delta\big(\Lu,\Lu\big)\to \mathcal{E}_\FF(\Lu)$ is the projection map.
\tpar In particular, let $f\in \FF T^\Delta\big(G,\Lu\big)$ be of the form
$f=\Ind_N^{G\times \Lu}\Inf_{\sur{N}}^NE$, 
where
\begin{itemize}
\item $N=N_{ G\times \Lu}\big(\Delta(P,\gamma,R)\big)$, for some $R\leq L$ and $\gamma:R\stackrel{\cong}{\to} P\leq G$,
\item $\sur{N}=N/\Delta(P,\gamma,R)$,
\item $E$ is a projective $k\sur{N}$-module.
\end{itemize}
If there exists some $v\in V$ such that $f\otimes v\notin \mathcal{R}$, then there exists $\varphi\in \FF T^\Delta(\Lu,G)$ such that $\pi(\varphi\circ f)\neq 0$. In particular, the Brauer quotient $(\varphi\circ f)[\Delta(L,\theta,L)]$ has to be non zero for some $\theta\in\Aut(L)$, which forces $R=L$.
\tpar Moreover if $f\in \FF T^\Delta\big(G,\Lu\big)$ can be factorized through a group of order strictly smaller than the order of $\Lu$, then $f\otimes v\in\mathcal{R}$ for any $v\in V$. So if there exists $v\in V$ with $f\otimes v\notin \mathcal{R}$, then in particular $p_2(N)=\Lu$. 
\begin{rem}{Remark}\label{conjugate} This condition $p_2(N)=\Lu$ means that for any $x\in\Lu$, there exists $g\in N_G(P)$ such that $^g\gamma(l)=\gamma({^xl})$ for all $l\in L$. Equivalently, there exists $s\in N_G(P)$ such that ${^s\gamma}(l)=\gamma(^ul)$ for all $l\in L$: Suppose indeed that such an element $s$ exists. Any element $x\in L$ is equal to $l_0u^\alpha$ for some $l_0\in L$ and $\alpha\in\NN$. Then
$$\gamma\circ i_x=\gamma\circ i_{l_0}\circ i_{u^\alpha}=i_{\gamma(l_0)}\circ\gamma\circ (i_u)^\alpha=i_{\gamma(l_0)}\circ (i_{s})^\alpha\circ\gamma=i_{\gamma(l_0)s^\alpha}\circ\gamma,$$
that is $i_g\circ \gamma=\gamma\circ i_x$, for $g=\gamma(l_0)s^\alpha$. \par
In other words, saying that $p_2(N)=\Lu$ amounts to saying that $(P,\gamma)$ belongs to the set
$$\index{\ordidx{PA}$\CP(G,L,u)$}
\CP(G,L,u):=\big\{(P,\gamma)\mid \gamma:L\stackrel{\cong}{\to}P\leq G,\;\exists s\in N_G(P),\;i_s\circ \gamma=\gamma\circ i_u\big\}.$$ 
\end{rem}
For $\varphi\in \Aut\big(\Lu\big)$, we have that
\begin{align*}
N_{G\times \Lu}\big(\Delta(P,\gamma\varphi,L)\big)&=\big\{(a,b)\mid (a,\varphi(b))\in N_{G\times \Lu}\big(\Delta(P,\gamma,L)\big)\big\}\\
&=(1\times\varphi^{-1})N_{G\times \Lu}\big(\Delta(P,\gamma,L)\big)\end{align*}
so the set $\CP(G,L,u)$ is a $\Big(G,\Aut\big(\Lu\big)\Big)$-biset by
$$\forall g\in G,\,\forall \varphi\in \Aut(\Lu),\;g\cdot (P,\gamma)\cdot\varphi=({^gP},i_g\gamma\varphi),$$
\tpar Let $\CT(G,L,u)$ \index{\ordidx{TV}$\CT(G,L,u)$} denote the set of triples $(P,\gamma,E)$, where $(P,\gamma)\in\CP(G,L,u)$ and $E$ is an indecomposable projective $kN\big(\Delta(P,\gamma,L)\big)/\Delta(P,\gamma,L)$-module. The set $\CT(G,L,u)$ is also a $\Big(G,\Aut\big(\Lu\big)\Big)$-biset by
$$\forall g\in G,\,\forall \varphi\in \Aut(\Lu),\;g\cdot(P,\gamma,E)\cdot\varphi=({^gP},i_g\gamma \varphi,{^gE}^\varphi)$$
where $^gE^\varphi$ is the $kN\big(\Delta({^gP},i_g\gamma\varphi,L)\big)/\Delta({^gP},i_g\gamma\varphi,L)$-module obtained from $E$ via the group isomorphism
$$(a,b)\Delta({^gP},i_g\gamma\varphi,L)\mapsto \big(a^g,\varphi(b)\big)\Delta(P,\gamma,L)$$
from $kN\big(\Delta({^gP},i_g\gamma\varphi,L)\big)/\Delta({^gP},i_g\gamma\varphi,L)$ to $kN\big(\Delta(P,\gamma,L)\big)/\Delta(P,\gamma,L)$.\mpn
For $(P,\gamma,E)\in\CT(G,L,u)$, set 
$$T(P,\gamma,E)\index{\ordidx{TT}$T(P,\gamma,E)$}=\Ind_N^{G\times\Lu}\Inf_{\sur{N}}^NE,$$
where $N=N_{G\times \Lu}\big(\Delta(P,\gamma,L)\big)$ and $\sur{N}=N/\Delta(P,\gamma,L)$. Then $T(P,\gamma,E)$ is a diagonal $p$-permutation $(kG,k\Lu)$-bimodule, and we abuse notation writing $T(P,\gamma,E)\in \FF T^\Delta(G,\Lu)$. We observe that for $(g,\varphi)\in G\times\Aut(\Lu)$, we have
$$T({^gP},i_g\gamma\varphi,{^gE^\varphi})\cong T(P,\gamma,E)\otimes_{\Lu}k\big(\Lu\big)_\varphi,$$
where $k\big(\Lu\big)_\varphi$ is the bimodule $k\big(\Lu\big)$ twisted by $\varphi$ on the right, that is $x\cdot m\cdot y\;\hbox{(in $k\big(\Lu\big)_\varphi$)}\,=xm\varphi(y)\;\hbox{(in $k\big(\Lu\big)$)}$.
\tpar Now a lemma:
\begin{enonce}{Lemma}
Let $J$ be a finite group, let $K\normal J$ such that $J/K$ is a cyclic $p'$-group. Let $E$ be an indecomposable projective $kJ$-module, let $V$ be an indecomposable summand of $\Res_K^JE$, and $H$ be the inertia group of $V$ in $J$. Then $V$ extends to an indecomposable projective $kH$-module $F$, and there exists a group homomorphism $\lambda :H/K\to k^\times$ such that 
$$E\cong\Ind_H^J(F\otimes_k \Inf_{H/K}^Hk_\lambda).$$
\end{enonce}
\pf Use Theorem 3.13.2 in~\cite{benson1}, and Theorem 4.1~of~\cite{functorialequivalence}.\findemo 
\tpar \label{Pim}Let $(P,\gamma,E)\in\CT(G,L,u)$. Set $N_{P,\gamma}\index{\ordidx{NA}$N_{P,\gamma}$}=N_{G\times \Lu}\big(\Delta(P,\gamma,L)\big)$, and $\sur{N}_{P,\gamma}\index{\ordidx{NB}$\sur{N}_{P,\gamma}$}=N_{P,\gamma}/\Delta(P,\gamma,L)$, so $E$ is an indecomposable projective $k\sur{N}_{P,\gamma}$-module. There is an exact sequence of groups
\begin{moneq}\label{suite exacte}
1\to C_G(P)\to \sur{N}_{P,\gamma}\to \langle u\rangle\to 1,
\end{moneq}
so by the previous lemma, if $V$ is an indecomposable summand of $\Res_{C_G(P)}^{\sur{N}_{P,\gamma}}E$, and $H$ its stabilizer in $\sur{N}_{P,\gamma}$, there exists an indecomposable projective $kH$-module $F$ such that  
$$E\cong\Ind_{H}^{\sur{N}_{P,\gamma}}F.$$
From this follows that
\begin{align*}
T(P,\gamma,E)&=\Ind_{N_{P,\gamma}}^{G\times \Lu}\Inf_{\sur{N}_{P,\gamma}}^{N_{P,\gamma}}E\\
&\cong\Ind_{N_{P,\gamma}}^{G\times \Lu}\Inf_{\sur{N}_{P,\gamma}}^{N_{P,\gamma}}\Ind_{H}^{\sur{N}_{P,\gamma}}F\\
&\cong\Ind_{N_{P,\gamma}}^{G\times \Lu}\Ind_{\hat{H}}^{N_{P,\gamma}}\Inf_H^{\hat{H}}F\cong\Ind_{\hat{H}}^{G\times \Lu}\Inf_H^{\hat{H}}F,
\end{align*}
where $\hat{H}$ is the inverse image in $N_{P,\gamma}$ of $H\leq \sur{N}_{P,\gamma}$ by the projection map $N_{P,\gamma}\to \sur{N}_{P,\gamma}$. Now $T(P,\gamma,E)$ factors through the second projection of $\hat{H}\leq G\times \Lu$, so we can assume that this projection is the whole of $\Lu$, i.e. equivalently that $H=\sur{N}_{P,\gamma}$. In other words, by Theorem 4.1 of~\cite{functorialequivalence} already quoted above, we can assume that $\Res_{C_G(P)}^{\sur{N}_{P,\gamma}}E$ is indecomposable. We denote by $\mathrm{Pim}^\sharp(k\sur{N}_{P,\gamma})$ \index{\ordidx{PP}$\mathrm{Pim}^\sharp(k\sur{N}_{P,\gamma})$} the set of isomorphism classes of such indecomposable projective $k\sur{N}_{P,\gamma}$-modules, and by $\CT^\sharp(G,L,u)$ \index{\ordidx{TV}$\CT^\sharp(G,L,u)$} the subset of $\CT(G,L,u)$ consisting of triples $(P,\gamma,E)$ such that $E\in\mathrm{Pim}^\sharp(k\sur{N}_{P,\gamma})$.
\tpar It follows from the above remarks that $S_{\Lu,V}(G)$ is generated by the images of the elements $T(P,\gamma,E)\otimes v$, where 
\begin{itemize}
\item $(P,\gamma)$ runs through a set $[\CP(G,L,u)]$ \index{\ordidx{PP}$[\CP(G,L,u)]$} of representatives of orbits of $\big(G\times\Aut(\Lu)\big)$ on $\CP(G,L,u)$, 
\item $E\in\mathrm{Pim}^\sharp(k\sur{N}_{P,\gamma})$, 
\item $v\in V$. 
\end{itemize}
In other words, we have a surjective map
$$\mathop{\bigoplus}_{\substack{(P,\gamma)\in [\CP(G,L,u)]\\E\in\mathrm{Pim}^\sharp(k\sur{N}_{P,\gamma})}}T(P,\gamma,E)\otimes V\longrightarrow S_{\Lu,V}(G)$$
sending $T(P,\gamma,E)\otimes v\in \F T^\Delta\big(G,\Lu\big)\otimes V$ to its image in $S_{\Lu,V}(G)$. The kernel of this map is the set of sums
$$\sum_{\substack{(P,\gamma)\in [\CP(G,L,u)]\\E\in\mathrm{Pim}^\sharp(k\sur{N}_{P,\gamma})}}T(P,\gamma,E)\otimes v_{P,\gamma,E}$$
where $v_{P,\gamma,E}\in V$, such that for any $(Q,\delta)$ in $[\CP(G,L,u)]$ and any indecomposable projective $k\sur{N}_{Q,\delta}$-module $F$, or equivalently for any $F$ in $\mathrm{Pim}^\sharp(k\sur{N}_{Q,\delta})$
\begin{moneq}\label{relation}
\sum_{\substack{(P,\gamma)\in [\CP(G,L,u)]\\E\in\mathrm{Pim}^\sharp(k\sur{N}_{P,\gamma})}}\pi\big(T^o(Q,\delta,F)\otimes_{kG}T(P,\gamma,E)\big)\cdot v_{P,\gamma,E}=0,
\end{moneq}
{\flushleft where} $T^o(Q,\delta,F)$ \index{\ordidx{TT}$T^o(Q,\delta,F)$} is the $(k\Lu,kG)$-bimodule ``opposite'' of the $(kG,k\Lu)$-bimodule $T(Q,\delta,F)$. In other words
$$T^o(Q,\delta,F)=\Ind_{N^o_{Q,\delta}}^{\Lu\times G}\Inf_{\sur{N}^o_{Q,\delta}}^{N^o_{Q,\delta}}F^o$$
where 
\begin{itemize}
\item $N^o_{Q,\delta}=N_{\Lu\times G}\big(\Delta(L,\delta^{-1},Q)\big)$,
\item $\sur{N}^o_{Q,\delta}=N^o_{Q,\delta}/\Delta(L,\delta^{-1},Q)$,
\item $F^o$ is the opposite module of $F$.
\end{itemize}

Now if $\pi\big(T^o(Q,\delta,F)\otimes_{kG}T(P,\gamma,E)\big)\neq 0$, there exists an automorphism $\theta$ of $\Lu$ such that
$$\big(T^o(Q,\delta,F)\otimes_{kG}T(P,\gamma,E)\big)[\Delta(L,\theta,L)]\neq 0.$$
Hence there exist $V\leq G$, $\alpha:V\stackrel{\cong}{\to} L$, and $\beta:L\stackrel{\cong}{\to}V$ such that $\theta=\alpha\beta$ and
$$T^o(Q,\delta,F)[\Delta(L,\alpha,V)]\neq 0\;\hbox{and}\;T(P,\gamma,E)[\Delta(V,\beta,L)]\neq 0.$$
So $\Delta(Q,\delta,L)$ is conjugate to $\Delta(V,\alpha^{-1},L)$ in $G\times \Lu$, and $\Delta(P,\gamma,L)$ is conjugate to $\Delta(V,\beta,L)$ in $G\times \Lu$. By Remark~\ref{conjugate}, this amounts to saying that there exist $g,h\in G$ such that $(V,\alpha^{-1})=(^gQ,i_g\delta)$ and $(V,\beta)=(^hP,i_h\gamma)$. Hence $P=V^h={^{h^{-1}g}}Q$ and $\theta=\delta^{-1}i_{g^{-1}}\gamma$. \par
In other words, setting $z=h^{-1}g$, we have $P={^zQ}$ and $\gamma=i_z\delta\theta$, that is $(P,\gamma)=z\cdot(Q,\delta)\cdot\theta$. As $(P,\gamma)$ and $(Q,\delta)$ are both in our set $[\CP(G,L,u)]$ of representatives of $G\times \Aut(\Lu)$-orbits on $\CP(G,L,u)$, this forces $P=Q$ and $\gamma=\delta$.\par
Now Equation~\ref{relation} reduces to the fact that for every $(Q,\delta)$ in $[\CP(G,L,u)]$ and any $F\in \PimS{Q,\delta}$
\begin{moneq}\label{relation1}
\sum_{E\in \PimS{Q,\delta}}\pi\big(T^o(Q,\delta,F)\otimes_{kG}T(Q,\delta,E)\big) \cdot v_{Q,\delta,E}=0.
\end{moneq}
It follows that
\begin{moneq}\label{relationseule}
S_{\Lu,V}(G)\cong\!\!\mathop{\bigoplus}_{(Q,\delta)\in[\CP(G,L,u)]}\!\!\Big(\big(\mathop{\oplus}_{E\in\PimS{Q,\delta}} T(Q,\delta,E)\otimes V\big)/\mathcal{R}_{Q,\delta}\Big)
\end{moneq}
where the relations $\mathcal{R}_{Q,\delta}$ are given by (\ref{relation1}) for every $F\in\Pim{Q,\delta}$.
\par
Now we set
$$G_{Q,\delta}\index{\ordidx{GA}$G_{Q,\delta}$}=\{g\in G\mid\exists x\in \Lu,\;(g,x)\in N_{G\times \Lu}\big(\Delta(Q,\delta,L)\big)\}.$$
With this notation, the $\big(\Lu,\Lu\big)$-bimodule $M=T^o(Q,\delta,F)\otimes_{kG}T(Q,\delta,E)$ is isomorphic to
\begin{align}\refstepcounter{monequation}
M&\cong\mathop{\bigoplus}_{g\in G_{Q,\delta}\dom G/G_{Q,\delta}} \Ind_{N^o_{Q,\delta}*{^{(g,1)}N_{Q,\delta}}}^{\Lu\times \Lu}\big(\Inf_{\sur{N}^o_{Q,\delta}}^{N^o_{Q,\delta}}F^o\mathop{\otimes}_{kC_G(Q,{^gQ})}{^{(g,1)}\Inf}_{\sur{N}_{Q,\delta}}^{N_{Q,\delta}}E\big)\nonumber\\
&=\mathop{\bigoplus}_{g\in G_{Q,\delta}\dom G/G_{Q,\delta}} \Ind_{N^o_{Q,\delta}*{N_{^gQ,i_g\delta}}}^{\Lu\times \Lu}\big(\Inf_{\sur{N}^o_{Q,\delta}}^{N^o_{Q,\delta}}F^o\mathop{\otimes}_{kC_G(Q,{^gQ})}\Inf_{\sur{N}_{^gQ,i_g\delta}}^{N_{^gQ,i_g\delta}}{^{(g,1)}E}\big).\label{eqn}
\end{align}
\tpar Now another lemma:
\begin{enonce}{Lemma}
Let $G,H,K$ be finite groups, let $Z\leq X\leq G\times H$  and $T\leq Y\leq H\times K$ be subgroups. Set $D=k_2(X)\cap k_1(Y)$. Then $X/Z\mathop{\times}_D\limits Y/T$ is a $X*Y$-set, and for $(u,v)Z\in X/Z$ and $(w,r)T\in Y/T$, the stabilizer of $(u,v)Z\mathop{\times}_D\limits (w,r)T$ in $X*Y$ is equal to $^{(u,v)}Z*{^{(w,r)}T}$.
\end{enonce}
\pf This is straightforward.\findemo
\tpar \label{sec}Let $M_g$ denote the term of the direct sum~(\ref{eqn}) indexed by $g\in G$, that is
$$M_g=\Ind_{N^o_{Q,\delta}*{N_{^gQ,i_g\delta}}}^{\Lu\times \Lu}\big(\Inf_{\sur{N}^o_{Q,\delta}}^{N^o_{Q,\delta}}F^o\mathop{\otimes}_{kC_G(Q,{^gQ})}\Inf_{\sur{N}_{^gQ,i_g\delta}}^{N_{^gQ,i_g\delta}}{^{(g,1)}E}\big).$$
To apply the previous lemma, we set
$$X=N^o_{Q,\delta},\;Y=N_{^gQ,i_g\delta},\;Z=\Delta^o(Q,\delta,L),\;T=\Delta(^gQ,i_g\delta,L).$$
Then $k_2(X)=C_G(Q)$ and $k_1(Y)=C_G(^gQ)$, thus $D=C_G(Q,^gQ)$. Moreover $Z\normal X$ and $T\normal Y$\par
Now since $E$ is a direct summand of $k\sur{N}_{Q,\delta}$, it follows that $\Inf_{\sur{N}_{^gQ,i_g\delta}}^{N_{^gQ,i_g\delta}}{^{(g,1)}E}$ is a direct summand of $kY/T$. Similarly $\Inf_{\sur{N}^o_{Q,\delta}}^{N^o_{Q,\delta}}F^o$ is a direct summand of $kX/Z$. Hence $\Inf_{\sur{N}^o_{Q,\delta}}^{N^o_{Q,\delta}}F^o\mathop{\otimes}_{kC_G(Q,{^gQ})}\Inf_{\sur{N}_{^gQ,i_g\delta}}^{N_{^gQ,i_g\delta}}{^{(g,1)}E}$ is a direct summand of
$$kX/Z\otimes_{kC_G(Q,{^gQ})}kY/T\cong k(X/Z)\mathop{\times}_{D}(Y/T).$$
By the previous lemma, for $(u,v)Z\in X/Z$ and $(w,r)T\in Y/T$, the stabilizer of $(u,v)Z\mathop{\times}_D\limits (w,r)T$ in $X*Y$ is equal to $^{(u,v)}Z*{^{(w,r)}T}=Z*T$. Hence the vertices of the indecomposable direct summands of $kX/Z\otimes_{kC_G(Q,{^gQ})}kY/T$ are subgroups of $Z*T$.\par
It follows that if $\pi(M_g)\neq 0$, then there exists $\theta\in\Aut(\Lu)$ such that $\Delta(L,\theta,L)$ is conjugate in $\Lu\times \Lu$ to a subgroup of $Z*T$. Up to changing $\theta$ by some inner automorphism of $\Lu$, we can assume that $\Delta(L,\theta,L)\leq Z*T$. But
\begin{align*}
Z*T&=\Delta^o(Q,\delta,L)*\Delta({^gQ},i_g\delta,L)\\
&=\left\{(a,b)\in \Lu\times \Lu\mid \exists c\in G,\;\Big\{\!\!\begin{array}{l}(c,a)\in \Delta(Q,\delta,L)\\(c,b)\in\Delta({^gQ},i_g\delta,L)\end{array}\right\}\\
&=\big\{\big(\delta^{-1}(c),\delta^{-1}(c^g)\big)\mid c\in Q\cap{^gQ}\big\}.
\end{align*}
Then $Z*T$ contains $\Delta(L,\theta,L)$ if and only if $Q={^qQ}$ and $\theta(l)=\delta^{-1}({^g\delta(l)})$ for any $l\in L$. In other words $g\in N_G(Q)$ and $\delta^{-1}i_g\delta$ is the restriction of $\theta$ to~$L$.\par
For $g\in N_G(Q)$, we set ${i_g}^\delta=\delta^{-1}i_g\delta\in\Aut(L)$ and 
$$\widehat{G}_{Q,\delta}\index{\ordidx{GC}$\widehat{G}_{Q,\delta}$}=\{g\in N_G(Q)\mid \exists \theta\in\Aut(\Lu), \;{i_g}^\delta=\theta_{|L}\}.$$ 
With this notation, we see that $\pi(M_g)=0$ unless $g\in \widehat{G}_{Q,\delta}$. \par
We also observe that $G_{Q,\delta}$ is a normal subgroup of $\widehat{G}_{Q,\delta}$, and we set 
$$\sur{G}_{Q,\delta}\index{\ordidx{GB}$\sur{G}_{Q,\delta}$}=\widehat{G}_{Q,\delta}/G_{Q,\delta}.$$
We observe moreover that for $g\in \widehat{G}_{G,\delta}$, there is a unique $\theta\in\Aut(\Lu)$ such that $\theta_{|L}=\dgd$, up to inner automorphism, thanks to Lemma~\ref{inner}.
\tpar \label{sec2}For $g\in\widehat{G}_{Q,\delta}$, we will denote by $\theta_g\in\Out(\Lu)$ \index{\ordidx{TW}$\theta_g$}the unique outer automorphism such that $(\theta_g)_{|L}=\dgd$, and by $\hat{\theta}_g$ \index{\ordidx{TX}$\hat{\theta}_g$} a representative of $\theta_g$ in $\Aut(\Lu)$. The map $g\in\widehat{G}_{Q,\delta}\mapsto \theta_g\in\Out(\Lu)$ is a group homomorphism, with kernel $G_{Q,\delta}$. In other words $\sur{G}_{Q,\delta}=\widehat{G}_{Q,\delta}/G_{Q,\delta}$ embeds in $\Out(\Lu)$.\par
Now let $(a,b)\in N_{Q,\delta}$, i.e. $(a,b)\in G\times \Lu$, and $^a\delta(l)=\delta(^bl)$ for all $l\in L$. If $g\in \widehat{G}_{Q,\delta}$ we claim that $\big({^ga},\hat{\theta}_g(b)\big)$ also lies in $N_{Q,\delta}$. Indeed for $l\in L$, setting $l'=\hat{\theta}_g^{-1}(l)$, we have $\delta(l)=i_g\delta(l')$, so
\begin{align*}
^{(^ga)}\delta(l)&={^{(^ga)}i_g}\delta(l)=ga\delta(l')a^{-1}g^{-1}\\
&=g\delta(^bl')g^{-1}=i_g\delta(^bl'),
\end{align*}
whereas
\begin{align*}
\delta({^{\hat{\theta}_g(b)}l})&=\delta\big({^{\hat{\theta}_g(b)}\hat{\theta}_g(l')}\big)\\
&=\delta\big(\hat{\theta}_g(^bl')\big)\\
&=\delta\delta^{-1}i_g\delta(^bl')=i_g\delta(^bl'),
\end{align*}
proving the claim.\par
In other words, the map $\Phi_g:\index{\ordidx{PQ}$\Phi_g$}(a,b)\mapsto \big({^ga},\hat{\theta}_g(b)\big)$ is an automorphism of $N_{Q,\delta}$. Moreover, it sends $\Delta(Q,\delta,L)$ to itself. Indeed, for $l\in L$, we have
$$\big({^g\delta}(l),\hat{\theta}_g(l)\big)=\big(i_g\delta(l),\delta^{-1}i_g\delta(l)\big)=\Big(\delta\big(\theta_g(l)\big),\theta_g(l)\Big).$$
It follows that $\Phi_g$ induces an automorphism $\bar{\Phi}_g$ of $\sur{N}_{Q,\delta}$. But there is a little more: Let $\ssur{N}_{Q,\delta}$ \index{\ordidx{NC}$\ssur{N}_{Q,\delta}$} denote the quotient $\sur{N}_{Q,\delta}/Z(Q)$. Then by the above lemma, for $g,h \in\widehat{G}_{Q,\delta}$, the composition $\Phi_g\Phi_h$ is equal to $\Phi_{gh}$ modulo an inner automorphism induced by an element of $Z(L)$. It follows that $\Phi_g$ induces an automorphism $\Phi_g^\flat$ of $\ssur{N}_{Q,\delta}$, and that $\Phi_g^\flat\Phi_h^\flat =\Phi_{gh}^\flat$. In other words $\widehat{G}_{Q,\delta}$ acts on $\ssur{N}_{Q,\delta}$.
 
\tpar Assuming now that $g\in \widehat{G}_{Q,\delta}$, we have
$$M_g=\Ind_{N^o_{Q,\delta}*{N_{Q,i_g\delta}}}^{\Lu\times \Lu}\big(\Inf_{\sur{N}^o_{Q,\delta}}^{N^o_{Q,\delta}}\,F^o\mathop{\otimes}_{kC_G(Q)}\Inf_{\sur{N}_{Q,i_g\delta}}^{N_{Q,i_g\delta}}\,{^{(g,1)}E}\big).$$
Moreover
\begin{align*}
N^o_{Q,\delta}*{N_{Q,i_g\delta}}&=\{(a,b)\in \Lu\times \Lu\mid \exists c\in G, (c,a)\in N_{Q,\delta}, (c,b)\in N_{Q,i_g\delta}\}\\
&=\left\{(a,b)\in \Lu\times \Lu\mid \exists c\in G,\,\forall l\in L,\Big\{\begin{array}{l}^c\delta(l)=\delta(^al)\\^{cg}\delta(l)={^g\delta(^bl)}\end{array}\right\}\\
&\leq \{(a,b))\in \Lu\times \Lu\mid \forall l\in L,\,{^a\big(\dgd(l)\big)}=\dgd(^bl)\}.
\end{align*}
The last group is the normalizer of $\Delta(L,\dgd,L)$ in $\Lu\times \Lu$. Conversely, if $(a,b)\in N_{\Lu\times \Lu}(\Delta(L,\dgd,L))$, since $(Q,\delta)\in \CP(G,L,u)$, there exists $c\in N_G(Q)$ such that $(c,a)\in N_{Q,\delta}$. In other words, we have
$$\forall l\in L,\;{^a\big(\dgd(l)\big)}=\dgd(^bl)\;\;\hbox{and}\;\;^c\delta(l)=\delta(^al).$$
It follows that
$$\forall l\in L,\;i_g\delta(^bl)=\delta\big(^a\dgd(l)\big)={^c\delta}\dgd(l)={^ci_g}\delta(l),$$
that is $(c,b)\in N_{Q,i_g\delta}$. Thus $(a,b)\in N^o_{Q,\delta}*{N_{Q,i_g\delta}}$, and this gives
$$N^o_{Q,\delta}*{N_{Q,i_g\delta}}=N_{\Lu\times \Lu}\big(\Delta(L,\dgd,L)\big).$$
To simplify the notation, we denote this group by $N_{L,u}(\dgd)$.\par
We also observe that $\Delta(L,\dgd,L)=\Delta^o(Q,\delta,L)*\Delta(Q,i_g\delta,L)$. We denote this group by $\Delta_{L,u}(\dgd)$, and we set 
$$\sur{N}_{L,u}(\dgd)=N_{L,u}(\dgd)/\Delta_{L,u}(\dgd).$$
Since $(L,u)$ is a $D^\Delta$-pair, this group fits into a short exact sequence of groups
\begin{moneq}\label{ses}
1\to Z(L)\to \sur{N}_{L,u}(\dgd)\to \langle u\rangle\to 1.
\end{moneq}
This sequence splits by the map sending $v\in\langle u\rangle$ to $\big(\hat{\theta}_g(v),v\big)\in\sur{N}_{L,u}(\dgd)$.
\tpar From the above discussion follows that
$$\pi(M)=\sum_{g\in \sur{G}_{Q,\delta}}\pi\left(\Ind_{N_{L,u}(\dgd)}^{\Lu\times \Lu}\big(\Inf_{\sur{N}^o_{Q,\delta}}^{N^o_{Q,\delta}}F^o\mathop{\otimes}_{kC_G(Q)}\Inf_{\sur{N}_{Q,i_g\delta}}^{N_{Q,i_g\delta}}{^{(g,1)}E}\big)\right).$$
The subgroup $\Delta_{L,u}(\dgd)$ of $N_{L,u}(\dgd)$ acts trivially on the module 
$$T(F,g,E)=\Inf_{\sur{N}^o_{Q,\delta}}^{N^o_{Q,\delta}}F^o\mathop{\otimes}_{kC_G(Q)}\Inf_{\sur{N}_{Q,i_g\delta}}^{N_{Q,i_g\delta}}{^{(g,1)}E},$$
so $T(F,g,E)$ is inflated from a $k\sur{N}_{L,u}(\dgd)$-module
$$T(F,g,E)=\Inf_{\sur{N}_{L,u}(\dgd)}^{N_{L,u}(\dgd)}\sur{T}(F,g,E),$$
where $\sur{T}(F,g,E)=F^o\mathop{\otimes}_{kC_G(Q)}\limits{^{(g,1)}E}$.\par
The discussion in Section~\ref{sec} above shows that the vertices of the indecomposable summands of $T(F,g,E)$ are equal to $\Delta_{L,u}(\dgd)$. In other words $\sur{T}(F,g,E)$ is a projective $\sur{N}_{L,u}(\dgd)$-module. In view of the split exact sequence~(\ref{ses}), we have that $\sur{N}_{L,u}(\dgd)\cong Z(L)\rtimes \langle u\rangle$, and
$$\sur{T}(F,g,E)\cong \mathop{\bigoplus}_{\lambda\in\uhat}
m_\lambda(F,g,E)\Ind_{\langle u\rangle}^{Z(L)\rtimes \langle u\rangle}k_\lambda,$$
for some multiplicities $m_\lambda(F,g,E)\in\mathbb{N}$, where $\uhat$ \index{\ordidx{UU}$\uhat$} is the set of group homomorphisms $\lambda:\langle u\rangle\to k^\times$. It is easy to check that
\begin{align}\label{mult}\refstepcounter{monequation}
m_\lambda(F,g,E)&=\dim_k\Hom_{Z(L)\rtimes \langle u\rangle}\big(\Inf_{\langle u\rangle}^{Z(L)\rtimes \langle u\rangle}k_\lambda, \sur{T}(F,g,E)\big)\nonumber\\
&=\dim_k\Hom_{\langle u\rangle}\Big(k_{\lambda},\big(\sur{T}(F,g,E)\big)^{Z(L)}\Big).
\end{align}
Now the image of the module $\Ind_{N_{L,u}(\dgd)}^{\Lu\times\Lu}\Inf_{\sur{N}_{L,u}(\dgd)}^{N_{L,u}(\dgd)}\Ind_{\langle u\rangle}^{Z(L)\rtimes \langle u\rangle}k_\lambda$ in the essential algebra $\CE_\FF(\Lu)$ is equal to $kG_{\hat{\theta}_g,\lambda^{-1}}$, where $\theta_g\in\Out(\Lu)$ and $\hat{\theta}_g\in\Aut(\Lu)$ are defined in Section~\ref{sec2}\footnote{Note that for $g\in\widehat{G}_{Q,\delta}$, the outer automorphism $\theta_g$ only depends on $gG_{Q,\delta}\in\sur{G}_{Q,\delta}$.}. So our relations $\mathcal{R}_{Q,\delta}$ of (\ref{relationseule}) become
\begin{moneq}\label{relationseule2}
\forall F\in\Pim{Q,\delta},\;\sum_{\substack{g\in \sur{G}_{Q,\delta}\\\lambda\in\uhat\\E\in\PimS{Q,\delta}}} m_\lambda(F,g,E)\,kG_{\theta_g,\lambda}\cdot v_E=0.
\end{moneq}
\tpar In order to understand these relations, we are going to change them a little bit, by replacing $F$ by its dual and the tensor product $-\otimes_{kC_G(Q)}-$ appearing in $T(F,g,E)$ by $\Hom_{kC_G(Q)}(-,-)$, in other words, by setting now
\begin{moneq}\label{hom}\sur{T}(F,g,E)\index{\ordidx{TT}$\sur{T}(F,g,E)$} =\Hom_{kC_G(Q)}(F,{^{(g,1)}E}).\end{moneq}
In particular, the action of $\sur{N}_{L,u}(\dgd)$ on $\sur{T}(F,g,E)$ is given as follows: for $(a,b)\in N_{L,u}(\dgd)$, choose $s\in G$ such that $(s,a)\in N_{Q,\delta}$. Then $(s^g,b)\in N_{Q,\delta}$, and for $\varphi\in \Hom_{kC_G(Q)}(F,{^{(g,1)}E})$, we have
\begin{moneq}\label{action}\forall f\in F,\;\big((a,b)\varphi\big)(f)=(s^g,b)\varphi\big((s,a)^{-1}f\big).
\end{moneq}
The action of $Z(L)$ on $\sur{T}(F,g,E)$ is simply given by multiplication
$$\forall z\in Z(L),\,\forall \varphi\in\Hom_{kC_G(Q)}(F,{^{(g,1)}E}),\,\forall f\in F,\;(z\varphi)(f)=\delta(z)^g\cdot \varphi(f),$$
where $\delta(z)\in Z(Q)$. It follows that 
\begin{align*}
\sur{T}(F,g,E)^{Z(L)}&=\Hom_{kC_G(Q)}\big(F,{^{(g,1)}E}^{Z(Q)}\big)\\
&\cong\Hom_{kC_G(Q)/Z(Q)}\big(F^{Z(Q)},{^{(g,1)}E}^{Z(Q)}\big),
\end{align*}
where the last isomorphism comes from the fact that since $F$ is projective,  the module of co-invariants $F_{Z(Q)}$ is isomorphic to the module of invariants $F^{Z(Q)}$. \par
Now in view of (\ref{mult}), we have to describe the action of $u$, that is the element $\big(\htheta_g(u),u\big)$ of $N_{L,u}(\dgd)$, on $\sur{T}(F,g,E)$. For this, we  use~(\ref{action}), and we choose $s\in G$ such that $\big(s,\hat{\theta}_g(u)\big)\in N_{Q,\delta}$. Now for $\varphi\in\Hom_{C_G(Q)}(F,{^{(g,1)}E})$, we have
$$\forall f\in F,\;(u\varphi)(f)=(s^g,u)\varphi\big((s,\hat{\theta}_g(u))^{-1}f\big).$$
For each $\lambda\in\uhat$, an element of $\Hom_{\langle u\rangle}\Big(k_\lambda,\Hom_{kC_G(Q)}\big(F^{Z(Q)},{^{(g,1)}E}^{Z(Q)}\big)\Big)$ is now determined by an element $\varphi\in \Hom_{kC_G(Q)}(F,{^{(g,1)}E})$ such that 
$$\forall f\in F,\;\lambda(u)\varphi(f)=(s^g,u)\varphi\Big(\big(s,\hat{\theta}_g(u)\big)^{-1}f\Big).$$
In other words
$$\forall f\in F,\;\varphi\Big(\big(s,\hat{\theta}_g(u)\big)f\Big)=\lambda(u)^{-1}(s^g,u)\varphi(f).$$
Since $(s^g,u)=\Phi_g^{-1}\big((s,\hat{\theta}_g(u)\big)$, we get that
$$\forall f\in F,\;\varphi\Big(\big(s,\hat{\theta}_g(u)\big)f\Big)=\lambda(u)^{-1}\Phi_g^{-1}\big((s,\hat{\theta}_g(u)\big)\varphi(f).$$
Now let $c\in C_G(Q)$. Then $(c,1)\in N_{Q,\delta}$, and 
$$\forall f\in F,\;\varphi\big((c,1)f\big)=(c^g,1)\varphi(f)=\lambda(1)^{-1}\Phi_g^{-1}\big((c,1)\big)\varphi(f).$$
Let $\sur{(a,b)}$ denote the image in $\sur{N}_{Q,\delta}$ of $(a,b)\in N_{Q,\delta}$. Now the element $\sur{\big(s,\hat{\theta}_g(u)\big)}$, together with the elements $\sur{(c,1)}$, for $c\in C_G(Q)$, generate the whole of $\sur{N}_{Q,\delta}$. It follows that for any $\sur{(a,b)}\in \sur{N}_{Q,\delta}$
$$\varphi\big(\sur{(a,b)}f\big)=\lambda\big(\hat{\theta}_g^{-1}(b)\big)^{-1}\bar{\Phi}_g^{-1}\big(\sur{(a,b)}\big)\varphi(f).$$
In other words $\varphi$ is a morphism of $k\sur{N}_{Q,\delta}$-modules from $F$ to the module $^{[g,\lambda]}E$\index{\ordidx{EE}$^{[g,\lambda]}E$}, equal to $E$ as a $k$-vector space, but with action of $\sur{(a,b)}\in\sur{N}_{Q,\delta}$ given by
\begin{moneq}\label{defaction}
\forall e\in E,\;\sur{(a,b)}\cdot e\;\hbox{(in $^{[g,\lambda]}E$)}:=\lambda\big(\hat{\theta}_g^{-1}(b)\big)^{-1}\bar{\Phi}_g^{-1}\big(\sur{(a,b)}\big)\cdot e\;\hbox{(in $E$)}.
\end{moneq}
It follows that
$$m_\lambda(F,g,E)=\dim_k\Hom_{\ssur{N}_{Q,\delta}}\big(F^{Z(Q)},{^{[g,\lambda]}E}^{Z(Q)}\big),$$
so our relations~(\ref{relationseule2}) become
\begin{moneq}\label{relationseule3}
\forall F\in\Pim{Q,\delta},\!\!\!\sum_{\substack{g\in \sur{G}_{Q,\delta}\\\lambda\in\uhat\\E\in\PimS{Q,\delta}}}\!\!\!\!\! \dim_k\Hom_{\ssur{N}_{Q,\delta}}\big(F^{Z(Q)},{^{[g,\lambda]}E}^{Z(Q)}\big)kG_{\theta_g,{\lambda}}\cdot v_E=0.
\end{moneq}
\begin{rem}{Remark} It should be noted that in~(\ref{defaction}), the coefficient $\lambda\big(\hat{\theta}_g^{-1}(b)\big)^{-1}$ does not depend on the choice of $\htheta_g\in\Aut(\Lu)$ in the class $\theta_g\in\Out(\Lu)$, as two different choices differ by an inner automorphism, so the corresponding values $\htheta_g(b)$ are conjugate. Hence we could write $\lambda\big(\theta_g^{-1}(b)\big)^{-1}$ instead of $\lambda\big(\hat{\theta}_g^{-1}(b)\big)^{-1}$.
\end{rem}
\tpar Suppose now that $g,h\in \widehat{G}_{Q,\delta}$, and $\lambda,\mu\in\uhat$. If $E$ is a projective $k\sur{N}_{Q,\delta}$-module, we claim that the $k\sur{N}_{Q,\delta}$-modules $^{[h,\mu]}\big({^{[g,\lambda]}E}\big)$ and $^{[hg,(\mu\circ\hat{\theta}_g)\cdot\lambda]}E$ are isomorphic. Indeed, there exists an element $w\in Z(L)$ such that $\htheta_{hg}=\htheta_h\htheta_gi_w$. Now $(1,w)\in N_{Q,\delta}$, and we can define $f:E\to E$ by $f(e)=\sur{(1,w^{-1})}\cdot e$. Then $f$ is clearly an automorphism of the $k$-vector space~$E$. \mpn
{\bf Claim:} {\em The map $f$ is an isomorphism of $k\sur{N}_{Q,\delta}$-modules from $^{[h,\mu]}\big({^{[g,\lambda]}E}\big)$ to $^{[hg,(\mu\circ\hat{\theta}_g)\cdot\lambda]}E$.}\mpn
\pf Indeed, for $e\in E$, let $^{[g,\lambda]}e$ denote the element $e$ of $^{[g,\lambda]}E$. Then for $(a,b)\in N_{Q,\delta}$
\begin{align*}
\sur{(a,b)}\cdot {^{[h,\mu]}\big(}{^{[g,\lambda]}e}\big)&=\mu\big(\hat{\theta}_h^{-1}(b)\big)\bar{\Phi}_h^{-1}\big(\sur{(a,b)}\big)\cdot {^{[g,\lambda]}e}\\
&=\mu\big(\hat{\theta}_h^{-1}(b)\big)\sur{(a^h,\hat{\theta}_h^{-1}(b))}\cdot {^{[g,\lambda]}e}\\
&=\mu\big(\hat{\theta}_h^{-1}(b)\big)\lambda\big(\hat{\theta}_g^{-1}\hat{\theta}_h^{-1}(b)\big)\sur{\big(a^{hg},\hat{\theta}_g^{-1}\htheta_h^{-1}(b)\big)}\cdot e.
\end{align*}
It follows that 
\begin{moneq}\label{fmachin}f\Big(\sur{(a,b)}\cdot {^{[h,\mu]}\big(}{^{[g,\lambda]}e}\big)\Big)=\mu\big(\hat{\theta}_h^{-1}(b)\big)\lambda\big(\hat{\theta}_g^{-1}\hat{\theta}_h^{-1}(b)\big)\sur{\big(a^{hg},w^{-1}\hat{\theta}_g^{-1}\htheta_h^{-1}(b)\big)}\cdot e.
\end{moneq}
On the other hand
\begin{align}\refstepcounter{monequation}
\sur{(a,b)}\cdot{^{[hg,(\mu\circ\hat{\theta}_g)\cdot\lambda]}f(e)}&=\big((\mu\circ\hat{\theta}_g)\cdot\lambda\big)\big(\hat{\theta}_{hg}^{-1}(b)\big)\sur{\big(a^{hg},\htheta_{hg}^{-1}(b)\big)}\cdot f(e)\nonumber\\
&=\mu\big(\htheta_g\htheta_{hg}^{-1}(b)\big)\lambda\big(\htheta_{hg}^{-1}(b)\big)\sur{\big(a^{hg},\htheta_{hg}^{-1}(b)w^{-1}\big)}\cdot e\label{machinf}
\end{align}
Since $\htheta_{hg}=\htheta_h\htheta_gi_w$, we have $i_w^{-1}\htheta_g^{-1}\htheta_h^{-1}=\htheta_{hg}^{-1}$, so $\htheta_{hg}^{-1}(b)w^{-1}=w^{-1}\htheta_g^{-1}\htheta_h^{-1}(b)$. Moreover 
$$\lambda\big(\htheta_{hg}^{-1}(b)\big)=\lambda\big(w^{-1}\hat{\theta}_g^{-1}\hat{\theta}_h^{-1}(b)w\big)=\lambda\big(\hat{\theta}_g^{-1}\hat{\theta}_h^{-1}(b)\big),$$
and 
\begin{align*}
\mu\big(\htheta_g\htheta_{hg}^{-1}(b)\big)&=\mu\Big(\htheta_g\big(w^{-1}\htheta_g^{-1}\htheta_h^{-1}(b)w\big)\Big)\\
&=\mu\Big(\htheta_g\big(w^{-1})\htheta_h^{-1}(b)\theta_g(w)\Big)=\mu\big(\htheta_h^{-1}(b)\big).
\end{align*}
Now it follows from (\ref{fmachin}) and (\ref{machinf}) that 
$$f\Big(\sur{(a,b)}\cdot {^{[h,\mu]}\big(}{^{[g,\lambda]}e}\big)\Big)=\sur{(a,b)}\cdot{^{[hg,(\mu\circ\hat{\theta}_g)\cdot\lambda]}f(e)},$$
proving the claim.\findemo
In addition to the above claim, we observe that if $g\in \widehat{G}_{Q,\delta}$ is actually in $G_{Q,\delta}$, then the $k\sur{N}_{Q,\delta}$-modules $E$ and $^{[g,1]}E$ are isomorphic: Indeed, saying that $g\in G_{Q,\delta}$ amounts to saying that $\htheta_g$ is an inner automorphism $i_x$ of $\Lu$, for some $x\in \Lu$ with $(g,x)\in N_{Q,\delta}$. Then~(\ref{defaction}) becomes 
\begin{align*}
\forall e\in E,\;\sur{(a,b)}\cdot e\;\hbox{(in $^{[g,1]}E$)}&=\bar{\Phi}_g^{-1}\big(\sur{(a,b)}\big)\cdot e\;\hbox{(in $E$)}\\
&=\sur{(a^g,b^x)}\cdot e.
\end{align*}
Then the map $e\in E\mapsto \sur{(g^{-1},x^{-1})}\cdot e\in {^{[g,1]}E}$ is an isomorphism of $k\sur{N}_{Q,\delta}$-modules.\medskip\par
We can now introduce the semidirect product
$$\semiG:=\big(\widehat{G}_{Q,\delta}/G_{Q,\delta}\big)\ltimes \uhat.$$
As a set $\semiG$ is the cartesian product $\sur{G}_{Q,\delta}\times \uhat$. For $g\in G_{Q,\delta}$ and $\lambda\in \uhat$, let $[g,\lambda]$ denote the pair $(gG_{Q,\delta},\lambda)$ in $\semiG$. The product in $\semiG$ is given as follows: For $g,h\in G_{Q,\delta}$ and $\lambda, \mu\in\uhat$, we have
$$[h,\mu][g,\lambda]=[hg,(\mu\circ\theta_g)\cdot\lambda].$$
The above discussion now shows that there is an action of $\semiG$ on the group $\mathrm{Proj}(k\sur{N}_{Q,\delta})$. The group $\semiG$ also acts on the set $\PimS{Q,\delta}$ introduced at Section~\ref{Pim}: Indeed the restriction of $^{[g,\lambda]}E$ to $C_G(Q)$ is isomorphic to the restriction of $^{[g,1]}E$, and $\widehat{G}_{Q,\delta}$ permutes the indecomposable $kC_G(Q)$-modules. 
\pagebreak[3]
\tpar Yet another (well known) lemma\footnote{We will not use Assertion 2 of this lemma here, but we state it for completeness.}:
\begin{enonce}{Lemma}
Let $H$ be a finite group, and $R$ be a normal $p$-subgroup of~$H$. \begin{enumerate}
\item The assignment $E\mapsto E^R$ induces a bijection between the set of isomorphism classes of indecomposable projective $kH$-modules and the set of isomorphism classes of indecomposable projective $k(H/R)$-modules.
\item Moreover if $R$ is central in $H$, then for any projective $kH$-modules $E$ and~$F$
$$\dim_k\Hom_{kH}(E,F)=|R|\dim_{k(H/R)}(E^R,F^R).$$
\end{enumerate}
\end{enonce}
\pf 1. Let $E$ be a projective $kH$-module, and $M$ be any $k(H/R)$-module. Then
$$\Hom_{kH}(E,\Inf_{H/R}^HM)\cong\Hom_{k(H/R)}(E_R,M)\cong \Hom_{k(H/R)}(E^R,M)$$
as $E_R\cong E^R$ if $E$ is projective. Now the functor $M\mapsto \Hom_{kH}(E,\Inf_{H/R}^HM)$ is exact, since inflation is exact, so $E^R$ is a projective $k(H/R)$-module. Moreover, the simple $kH$-modules are inflated from $H/R$, and the previous isomorphism shows that if $E$ is indecomposable, then $E^R$ has a unique simple $k(H/R)$-quotient, thus it is indecomposable.\mpn
2. If now $R$ is central  in $H$, then $R$ acts on $\Hom_{kH}(E,F)$ by left multiplication, that is $(r\varphi)(e)=r\cdot\varphi(e)$ for $r\in R$, $e\in E$, and $\varphi\in \Hom_{kH}(E,F)$. Moreover, this action is free, since if $E=F=kH$, then $\Hom_{kH}(E,F)\cong kH$ is free as a $kR$-module. Thus
\begin{align*}
\dim_k\Hom_{kH}(E,F)&=|R|\dim_k\big(\Hom_{kH}(E,F)\big)^R=|R|\dim_k\Hom_{kH}(E,F^R)\\
&=|R|\dim_k\Hom_{k(H/R)}(E_R,F^R)\\
&=|R|\dim_k\Hom_{k(H/R)}(E^R,F^R),
\end{align*}
as was to be shown.  \findemo
\tpar Let $\F R_k(\ssur{N}_{Q,\delta})\index{\ordidx{FR}$\F R_k(\ssur{N}_{Q,\delta})$} =\F\otimes_\Z R_k(\ssur{N}_{Q,\delta})$ denote the Grothendieck group of finite dimensional $k\ssur{N}_{Q,\delta}$-modules, extended by $\F$. For a projective $k\ssur{N}_{Q,\delta}$-module $X$, let $[X]$ denote its image in $\F R_k(\ssur{N}_{Q,\delta})$. The group $\semiG$ also acts on $R_k(\ssur{N}_{Q,\delta})$, by permutation of its base consisting of the isomorphism classes of simple modules.\par
We denote by $\Gamma_{Q,\delta}$\index{\ordidx{GH}$\Gamma_{Q,\delta}$} the image of the linear map 
$$\gamma^{L,u}_{Q,\delta}\index{\ordidx{GK}$\gamma^{L,u}_{Q,\delta}$} :\F\mathrm{Proj}^\sharp(k\sur{N}_{Q,\delta})\to \F R_k(\ssur{N}_{Q,\delta})$$
induced by $E\mapsto [E^{Z(Q)}]$. The map $\gamma_{Q,\delta}^{L,u}$ is a map of $\F (\semiG)$-modules, so its image $\Gamma_{Q,\delta}$ is also a $\F (\semiG)$-module.\par
The $\CE_\F(\Lu)$-module $V$ is also a $\F(\semiG)$-module, thanks to the (surjective) homomorphism of algebras sending $[g,\lambda]\in\F(\semiG)$ to the image of $kG_{\theta_g,\lambda}$ in $\CE_\F(\Lu)$.
Tensoring with $V$ gives a surjective map 
$$\gamma^{L,u}_{Q,\delta}\otimes \mathrm{Id}_V:\F\mathrm{Proj}^\sharp(k\sur{N}_{Q,\delta})\otimes_k V\to \Gamma_{Q,\delta}\otimes_k V,$$
of $\F(\semiG)$-modules, where the action of $\semiG$ on both tensor products is diagonal.\par
Now rephrasing~(\ref{relationseule}), we get a surjective map
$$\sigma=\mathop{\bigoplus}\sigma_{Q,\delta}: \mathop{\bigoplus}_{(Q,\delta)\in[\CP(G,l,u)]}\ProjS{Q,\delta}\otimes_k V\to S_{\Lu,V}(G),$$
where $\sigma_{Q,\delta}$ sends $E\otimes v\in\ProjS{Q,\delta}\otimes_k V$ to the image of $T(Q,\delta,E)\otimes v$ in $S_{\Lu,V}(G)$. \par
The kernel of this map is the direct sum for $(Q,\delta)\in[\CP(G,L,u)]$ of the kernels of its components $\sigma_{Q,\delta}$. By~(\ref{relationseule3}), the sum $\sum_{E\in\PimS{Q,\delta}}\limits E\otimes v_E$ is in the kernel of $\sigma_{Q,\delta}$ if and only if 
\begin{moneq}\label{kernel}
\forall F\in\Pim{Q,\delta},\!\!\!\!\sum_{\substack{[g,\lambda]\in\semiG\\E\in\PimS{Q,\delta}}}\!\!\!\! \dim_k\Hom_{\ssur{N}_{Q,\delta}}\big(F^{Z(Q)},{^{[g,\lambda]}E}^{Z(Q)}\big)[g,\lambda]\cdot v_E=0.
\end{moneq}
Now $R_k(\ssur{N}_{Q,\delta})$ has a basis consisting of the $[X_F]$ for $F\in\Pim{Q,\delta}$, where $X_F$ is the unique simple quotient of the $k\ssur{N}_{Q,\delta}$-module $F^{Z(Q)}$.
So (\ref{kernel}) is equivalent to
$$\sum_{F\in\Pim{Q,\delta}}[X_F]\otimes\sum_{\substack{[g,\lambda]\in\semiG\\E\in\PimS{Q,\delta}}}\dim_k\Hom_{\ssur{N}_{Q,\delta}}\big(F^{Z(Q)},{^{[g,\lambda]}E}^{Z(Q)}\big)[g,\lambda]\cdot v_E=0$$
in $\F R_k(\sur{N}_{Q,\delta})\otimes V$. This in turn is equivalent to
$$\sum_{F\in\Pim{Q,\delta}}\sum_{\substack{[g,\lambda]\in\semiG\\E\in\PimS{Q,\delta}}}\dim_k\Hom_{\ssur{N}_{Q,\delta}}\big(F^{Z(Q)},{^{[g,\lambda]}E}^{Z(Q)}\big)[X_F]\otimes [g,\lambda]\cdot v_E=0.$$
Now for any $[g,\lambda]\in \semiG$ and any $E\in \PimS{Q,\delta}$
$$\sum_{F\in\Pim{Q,\delta}}\dim_k\Hom_{\ssur{N}_{Q,\delta}}\big(F^{Z(Q)},{^{[g,\lambda]}E}^{Z(Q)})[X_F]=[{^{[g,\lambda]}E}^{Z(Q)}]={^{[g,\lambda]}\gamma}_{Q,\delta}^{L,u}(E).$$
It follows that $\sum_{E\in\PimS{Q,\delta}}\limits E\otimes v_E$ is in the kernel of $\sigma_{Q,\delta}$ if and only if
$$\sum_{E\in\PimS{Q,\delta}}\sum_{[g,\lambda]\in\semiG}{^{[g,\lambda]}\gamma}_{Q,\delta}^{L,u}(E)\otimes [g,\lambda]\cdot v_E=0$$
in $\F R_k(\sur{N}_{Q,\delta})\otimes V$. In other words $\sigma_{Q,\delta}$ has the same kernel as the map 
$$\sum_{E\in\PimS{Q,\delta}}E\otimes v_E\mapsto \sum_{[g,\lambda]\in\semiG}[g,\lambda]\big(\gamma_{Q,\delta}^{L,u}\otimes \mathrm{Id}\big)\big(\sum_{E\in\PimS{Q,\delta}}E\otimes v_E\big).$$
It follows that the image of $\sigma_{Q,\delta}$ is isomorphic to the image of the previous map, that is
$$\mathrm{Im}(\sigma_{Q,\delta})\cong\mathrm{Tr}_1^{\semiG}\big(\Gamma_{Q,\delta}\otimes V\big).$$
We finally get:
\begin{enonce}{Theorem} \label{isomorphism}There is an isomorphism of $k$-vector spaces
$$S_{\Lu,V}(G)\cong \mathop{\bigoplus}_{(Q,\delta)\in[\CP(G,L,u)]}\mathrm{Tr}_1^{\semiG}\big(\Gamma_{Q,\delta}\otimes V\big).$$
\end{enonce}
\tpar{\bf The simple diagonal $p$-permutation functors.} Now we use Theorem~\ref{isomorphism} to describe the simple functors. We assume that $\F$ is algebraically closed, of characteristic 0 or $p$. Recall that we have an isomorphism of algebras
$$\CE_\F\big(\Lu\big)\cong \F\sur{R}_k(\langle u\rangle)\rtimes \Out(\Lu).$$
In order to describe the simple $\CE_\F(\Lu)$-modules, we set $A=\F\sur{R}_k(\langle u\rangle)$ and $\Omega=\Out(\Lu)$, and we use the results of~\cite{montgomery-witherspoon}, in our specific situation, as in~Section 4.2 of~\cite{ducellier}. First, the simple $A$-modules are one dimensional, of the form $k_x$, where $x$ is a generator of $\langle u\rangle$: For such a generator $x$, we get an algebra homomorphism $e_x:A=\F\sur{R}_k(\langle u\rangle)\to \F$ defined by
$$\forall \lambda\in\langle u\rangle^\natural, \;e_x(\lambda)=\widetilde{\lambda}(x),$$
where $\widetilde{\lambda}:\langle u\rangle\to\F^\times$ lifts $\lambda$. This makes sense because if $\lambda$ is induced from a proper  subgroup of $\langle u\rangle$, then $\widetilde{\lambda}(x)=0$, as $x$ cannot be contained in a proper subgroup of $\langle u\rangle$. So $e_x$ extends to an algebra homomorphism $A=\F\sur{R}_k(\langle u\rangle)\to \F$, which in turn yields a one dimensional $A$-module $k_x$. We get all the simple $A$-modules in this way.\par
Now the stabilizer $\Omega_x$ of $k_x$ in $\Omega=\Out(\Lu)$ is the set of classes of $\gamma\in \Aut(\Lu)$ such that $\gamma(x)=x$. This does not depend on the generator $x$ of $\langle u\rangle$, and it is equal to $\Out(L,u)$. We note that $\Omega_x$ acts trivially on $A$, since it acts trivially on $\langle u\rangle$, so $A\rtimes \Omega_x=A\times\Omega_x$.\par
So any simple $\CE_\F\big(\Lu\big)$-module $V$ if of the form 
$$V=\Ind_{A\times \Omega_x}^{A\rtimes\Omega}(k_x\otimes W),$$
for some generator $x$ of $\langle u\rangle$ and some simple $\F\Omega_x$-module $W$. The action of $A\times\Omega_x$ on $k_x\otimes W$ is given by
$$\forall a\in A, \,\forall \gamma\in\Omega_x, \,\forall w\in W,\;(a,\gamma)\cdot (1\otimes w)=e_x(a)(1\otimes\gamma\cdot w).$$
Moreover, the isomorphism type of the simple $\CE_\F(\Lu)$-module $V$ is determined by the isomorphism type of the simple $\F\Out(L,u)$-module $W$, and the choice of the generator $x$ of $\langle u\rangle$, up to the action of $\Omega$, i.e. up to the action of the subgroup $\Aut(\Lu)^\sharp$ of $\Aut(\Lu)$ consisting of automorphisms which stabilize $\langle u\rangle$. \par
Then $(L,x)$ is a $D^\Delta$-pair, such that $\Out(L,u)=\Out(L,x)$, and choosing $x$ up to the action of $\Aut(\Lu)^\sharp$ amounts to choosing $(L,x)$ in a set of $D^\Delta$-pairs such that $L\langle x\rangle=\Lu$, up to isomorphism of $D^\Delta$-pairs. So up to changing $(L,u)$ to $(L,x)$, we can parametrize the simple functor $S_{\Lu,V}$ by the triple $(L,x,W)$ instead, that is, we can suppose $x=u$ in the previous calculations, and set $\SS_{L,u,W}\index{\ordidx{SS}$\SS_{L,u,W}$} =S_{\Lu,V}$.\par
Now we denote by $G_{Q,\delta,u}$ \index{\ordidx{GD}$G_{Q,\delta,u}$} the subgroup of $\widehat{G}_{Q,\delta}$ defined by

\begin{align*}
G_{Q,\delta,u}&=\{g\in N_G(Q)\mid \exists \theta\in \Aut(L,u), i_g^\delta=\theta_{|L}\}\\
&=\{g\in N_G(Q)\mid\theta_g\in\Aut(L,u)\},
\end{align*}
and we set
$$\sur{G}_{Q,\delta,u}\index{\ordidx{GD}$\sur{G}_{Q,\delta,u}$}=G_{Q,\delta,u}/G_{Q,\delta}.$$
With this notation, we have an isomorphism of $k\big(\sur{G}_{Q,\delta}\ltimes\langle  u\rangle^\natural\big)$-modules
$$\Gamma_{Q,\delta}\otimes V\cong\Ind_{\sur{G}_{Q,\delta,u}\times \langle u\rangle^\natural}^{\sur{G}_{Q,\delta}\ltimes\langle  u\rangle^\natural}\Big(\big(\Res\,\Gamma_{Q,\delta}\big)\otimes (k_u\otimes W)\Big),$$
where $\Res \,\Gamma_{Q,\delta}=\Res_{\sur{G}_{Q,\delta,u}\times \langle u\rangle^\natural}^{\sur{G}_{Q,\delta}\ltimes\langle  u\rangle^\natural}\Gamma_{Q,\delta}$, and the action of $\sur{G}_{Q,\delta,u}\times \langle u\rangle^\natural$ on $k_u\otimes W$ is given by
$$\forall g\in \sur{G}_{Q,\delta,u},\, \forall \lambda\in\langle u\rangle^\natural,\,\forall w\in W,\; (g,\lambda)\cdot (1\otimes w)=\lambda(u)\big(1\otimes \theta_g\cdot w\big).$$
It follows that
$$\mathrm{Tr}_1^{\sur{G}_{Q,\delta}\ltimes\langle  u\rangle^\natural}(\Gamma_{Q,\delta}\otimes V)\cong \mathrm{Tr}_1^{\sur{G}_{Q,\delta,u}\times\langle  u\rangle^\natural}\Big(\big(\Res\,\Gamma_{Q,\delta}\big)\otimes (k_u\otimes W)\Big).$$
Now for $\nu\in \Gamma_{Q,\delta}$ and $w\in W$, we have
\begin{align*}
\mathrm{Tr}_1^{\sur{G}_{Q,\delta,u}\times\langle  u\rangle^\natural}\big(\nu\otimes(1\otimes w)\big)&=\sum_{g\in \sur{G}_{Q,\delta,u}}\sum_{\lambda\in\langle u\rangle^\natural}(g,\lambda)\cdot\big(\nu\otimes(1\otimes w)\big)\\
&=\sum_{g\in \sur{G}_{Q,\delta,u}}\sum_{\lambda\in\langle u\rangle^\natural}\big((g,\lambda)\cdot \nu\big)\otimes \lambda(u)\big(1\otimes \theta_g\cdot w\big)\\
&=\sum_{g\in \sur{G}_{Q,\delta,u}}\big((g,1)\cdot\sum_{\lambda\in\langle u\rangle^\natural}\lambda(u)\nu_\lambda\big)\otimes \big(1\otimes \theta_g\cdot w\big).
\end{align*}
Now when $\nu$ runs through $\Gamma_{Q,\delta}$, the sum $\sum_{\lambda\in\langle u\rangle^\natural}\limits\lambda(u)\nu_\lambda$ runs through the subspace
$$\Gamma^0_{Q,\delta}\index{\ordidx{GK}$\Gamma^0_{Q,\delta}$}=\{m\in\Gamma_{Q,\delta}\mid \forall\lambda\in\langle u\rangle^\natural,\;m_\lambda=\widetilde{\lambda}(u)^{-1}m\}$$
of $\Gamma_{Q,\delta}$. This subspace is the image by the map $\gamma_{Q,\delta}^{L,u}$ of the corresponding subspace
$$\F \mathrm{Proj}^\sharp(k\sur{N}_{Q,\delta})^0\index{\ordidx{FP}$\F \mathrm{Proj}^\sharp(k\sur{N}_{Q,\delta})^0$}=\{n\in \F \mathrm{Proj}^\sharp(k\sur{N}_{Q,\delta})\mid\forall\lambda\in\langle u\rangle^\natural,\;n_\lambda=\widetilde{\lambda}(u)^{-1}n\}$$
of $\F \mathrm{Proj}^\sharp(k\sur{N}_{Q,\delta})$. This space has a basis consisting of the sums $\sum_{\lambda\in\langle u\rangle^\natural}\limits\lambda(u)E_\lambda$, where $E$ runs through a set $\Sigma$ of representatives of orbits of $\mathrm{Pim}^\sharp(k\sur{N}_{Q,\delta})$ under the action of $\langle u\rangle^\natural$. Since $\mathrm{Pim}^\sharp(k\sur{N}_{Q,\delta})$ is a $(\sur{G}_{Q,\delta,u},\langle u\rangle^\natural)$-biset, this set of representatives can moreover be chosen invariant by the action of $\sur{G}_{Q,\delta,u}$.\par
Now the restriction map $\Res_{C_G(Q)}^{\sur{N}_{Q,\delta}}$ induces a bijection from $\Sigma$ to the set $\mathrm{Pim}\big(kC_G(Q),u\big)$\index{\ordidx{PP}$\mathrm{Pim}\big(kC_G(Q),u\big)$} of isomorphism classes of $u$-invariant indecomposable projective $kC_G(Q)$-modules, and this bijection is $\sur{G}_{Q,\delta,u}$-equivariant. Moreover, taking fixed points by $Z(Q)$ gives a bijection from $\mathrm{Pim}\big(kC_G(Q),u\big)$ to the set $\mathrm{Pim}\big(kC_G(Q)/Z(Q),u\big)$ of isomorphism classes of $u$-invariant indecomposable projective $kC_G(Q)/Z(Q)$-modules.\par
Finally, we have proved the following:
\begin{enonce}{Theorem} \label{the simple functors}Let $\F$ be a field of characteristic 0 or $p$. 
\begin{enumerate}
\item The simple diagonal $p$-permutation functors over $\F$ are parametrized by triples $(L,u,W)$, where $(L,u)$ is a $D^\Delta$-pair and $W$ is a simple $\F\Out(L,u)$-module. 
\item The evaluation at a finite group $G$ of the simple functor $\SS_{L,u,W}$ parame\-trized by the triple $(L,u,W)$ is 
$$\SS_{L,u,W}(G)\cong \mathop{\bigoplus}_{(Q,\delta)\in[\CP(G,L,u)]}\mathrm{Tr}_1^{\sur{G}_{Q,\delta,u}}\Big(\F\Cart\big(kC_G(Q)/Z(Q),u\big)\otimes W\Big),$$
where $\F\Cart\big(kC_G(Q)/Z(Q),u\big)$ is the image of the map 
$$\F \mathrm{Pim}\big(kC_G(Q)/Z(Q),u\big)\to\F R_k\big(C_G(Q)/Z(Q)\big)$$ 
induced by the Cartan map.
\end{enumerate}
\vspace{-3ex}
\end{enonce}
\section{Examples}
\tpar {\bf The functor $\SS_{\un,1,\F}$.} We apply Theorem~\ref{the simple functors} to the case where $L=\un$, so $u=1$, and $W=\F$. For a finite group $G$, we get that $(Q,\delta)\in \CP(G,\un,1)$ if and only if $Q=\un$ and $\delta:L\to Q$ is the identity. Moreover $G_{Q,\delta,1}=G=G_{Q,\delta}$, so by Theorem~\ref{the simple functors}, we get that
$$ \SS_{\un,1,\F}(G)\cong\F\Cart(G),$$
where $\F\Cart(G)$\index{\ordidx{FC}$\F\Cart(G)$} is the image of the map $\F\Proj(kG)\to \F R_k(G)$ induced by the Cartan map $\sfc^G: \Proj(kG)\to R_k(G)$\index{\ordidx{CC}$\sfc^G$}. We now show that the previous isomorphism is quite natural.

For this, we observe that the assignments $G\mapsto \Proj(kG)$ and $G\mapsto R_k(G)$ are diagonal $p$-permutation functors: If $M$ is a diagonal $p$-permutation $(kH,kG)$-bimodule, then $M$ is left and right projective, so if $\Lambda$ is a projective $kG$-module, then $M\otimes_{kG}\Lambda$ is a projective $kH$-module. Similarly, the functor $M\otimes_{kG} -$ changes a short exact sequence of $kG$-modules into a short exact sequence of $kH$-modules.\par
Moreover, the Cartan maps $\sfc^G$ form a morphism of diagonal $p$-permuta\-tion functors
$$\sfc:\Proj(k-)\to R_k(-).$$
In particular the assignment $\F\Cart(-): G\mapsto \F\mathrm{Cart}(G)$ is a subfunctor of $\F R_k(-)$. 
\begin{enonce}{Lemma} \label{unique minimal}The functor $\F\Cart(-)$ is the unique minimal subfunctor of $\F R_k(-)$. It is isomorphic to the simple functor $\SS_{\un,1,\F}$. 
\end{enonce}
\pf Let $F$ be a subfunctor of $\F R_k(-)$. Then $F(\un)\leq \F R_k(\un)=\F$, so $F(\un)$ is either 0 or $\F$. Suppose first that $F(\un)=0$. Let $G$ be a finite group, and $u\in F(G)$. Then $u=\sum_{S\in\Irr_k(G)}\limits\lambda_SS$, where $\lambda_S\in\F$. Let $T\in \Irr_k(G)$, and $\mathsf{P}_T$ be its projective cover. Then $\sfP_T\in\F T^\Delta(\un,G)$, so $\sfP_T\otimes_{kG}u\in F(1)=0$. But for $S\in \Irr_k(G)$, we have $\sfP_T\otimes_{kG}S=0$ unless $S$ is isomorphic to the dual $T^\natural$ of $T$. It follows that $\lambda_{T^\natural}=0$ for any $T\in \Irr_k(G)$, so $u=0$. Hence $F=0$ if $F(\un)=0$.\par
Suppose now that $F(\un)=\F$, that is $F(\un)\ni k$. If $T\in \Irr_k(G)$, then $\sfP_T\in\F T^\Delta(G,\un)$, so $\sfP_T\otimes_{k}k\in F(G)$, for any $T\in \Irr_k(G)$. But $\sfP_T\otimes_{k}k\cong \sfP_T$ is the image of $\sfP_T$ by the Cartan map. It follows that $F(G)$ contains $\F\Cart(G)$, so $F\geq \F\Cart(-)$. Hence $\F\Cart(-)$ is the unique (non zero) minimal subfunctor of $\F R_k$. Since $\F\Cart(\un)\cong\F$, it follows that $\F\Cart(-)\cong \SS_{\un,1,\F}$, completing the proof.\endpf
\begin{rem}{Remark} When $\F$ has characteristic $p$, the functor $\F\Cart(-)$ is a proper subfunctor of $\F R_k(-)$, so Lemma~\ref{unique minimal} shows in particular that the category $\CF_{\F pp_k}^\Delta$ is {\em not} semisimple.
\end{rem}
\tpar In the case $\F$ has characteristic 0, the Cartan matrix has non zero determinant in $\F$, so the Cartan map $\F \sfc^G:\F\Proj(kG)\to \F R_k(G)$ is invertible. So we have isomorphisms of functors
$$\F\Proj(-)\cong \F R_k(-)\cong\SS_{\un,1,\F}$$
in this case. This is Theorem 5.20 in \cite{bouc-yilmaz}.
\tpar The other case we can consider is when $\F$ is a field of characteristic~$p$, and we assume that $\F=k$. 
We choose a $p$-modular system $(K,\CO,k)$, and we assume that $K$ is big enough for the group $G$. If $S$ is a simple $kG$-module, we denote by $\Phi_S:G_{p'}\to \CO$ the modular character of $\mathsf{P}_S$, where $G_{p'}$ is the set of $p$-regular elements of~$G$. If $v=\sum_{S\in \Irr_k(G)}\limits \omega_S\mathsf{P}_S$, where $\omega_S\in \CO$, is an element of $\CO\Proj(kG)$, we denote by $\Phi_v:\CO\Proj(kG)\to \CO$ the map $\sum_{S\in\Irr_k(G)}\limits \omega_S\Phi_S:G_{p'}\to \CO$, and we  call $\Phi_v$ the modular character of $v$. \par
Then for a simple $kG$-module $T$, the multiplicity of $S$ as a composition factor of $\mathsf{P}_T$ is equal to the Cartan coefficient
$$\sfc^G_{T,S}=\dim_k\Hom_{kG}(\mathsf{P}_T,\mathsf{P}_S)=\frac{1}{|G|}\sum_{x\in G_{p'}}\Phi_T(x)\Phi_S(x^{-1}).$$
In order to describe the image of the Cartan map $k\sfc^G$, we want to evaluate the image of this integer under the projection map $\rho:\CO\to k$. For this, we denote by $[G_{p'}]$ a set of representatives of conjugacy classes of $G_{p'}$, and we observe that in the field~$K$, we have
\begin{align}\refstepcounter{monequation}
\sfc^G_{S,T}&=\frac{1}{|G|}\sum_{x\in [G_{p'}]}\frac{|G|}{|C_G(x)|}\Phi_T(x)\Phi_S(x^{-1})\nonumber\\
&=\sum_{x\in [G_{p'}]}\frac{1}{|C_G(x)|}\frac{\Phi_T(x)}{|C_G(x)|_{p}}\frac{\Phi_S(x^{-1})}{|C_G(x)|_{p}}|C_G(x)|_{p}^2\nonumber\\
&=\sum_{x\in [G_{p'}]}\frac{\big(\Phi_T(x)/|C_G(x)|_{p}\big)\big(\Phi_S(x^{-1})/|C_G(x)|_{p}\big)}{|C_G(x)_{p'}|}|C_G(x)|_{p}.\label{in O}
\end{align}
But since $\Phi_S$ and $\Phi_T$ are characters of projective $kG$-modules (and since $C_G(x)=C_G(x^{-1}))$, the quotients $\Phi_T(x)/|C_G(x)|_{p}$ and $\Phi_S(x^{-1})/|C_G(x)|_{p}$ are in $\CO$, so 
$$\forall x\in[G_{p'}],\;\frac{\big(\Phi_T(x)/|C_G(x)|_{p}\big)\big(\Phi_S(x^{-1})/|C_G(x)|_{p}\big)}{|C_G(x)_{p'}|}\in\CO.$$
Now it follows from~\ref{in O} that
\begin{moneq}\label{rho}\rho(\sfc^G_{T,S})=\sum_{x\in [G_0]}\rho\left(\frac{\Phi_T(x)\Phi_S(x^{-1})}{|C_G(x)|}\right),
\end{moneq}
{\flushleft where} $[G_0]$ is a set of representatives of conjugacy classes of the set $G_0$ of elements {\em defect zero} of $G$, i.e. the set of $p'$-elements $x$ of $G$ such that $C_G(x)$ is a $p'$-group.
\begin{enonce}{Notation} \label{notation Gamma}For $x\in G_0$, we set
$$\Gamma_{G,\,x}=\sum_{S\in\Irr(kG)}\frac{\Phi_S(x^{-1})}{|C_G(x)|}S\in\CO R_k(G),$$
where $\Irr(kG)$ is a set of representatives of isomorphism classes of simple $kG$-modules. We also set
$$\gamma_{G,\,x}=\sum_{S\in\Irr(kG)}\rho\left(\frac{\Phi_S(x^{-1})}{|C_G(x)|}\right)S\in kR_k(G),$$
\end{enonce}
\begin{rem}{Remark}\label{class}
We note that $\Gamma_{G,\,x}$ and $\gamma_{G,\,x}$ only depend on the conjugacy class of $x$ in $G$, that is $\Gamma_{G,\,x}=\Gamma_{G,\,x^g}$ and $\gamma_{G,\,x}=\gamma_{G,\,x^g}$ for $g\in G$.
\end{rem}
{\flushleft By} Theorem 6.3.2 of~\cite{brauer-nesbitt} (see also Theorem 6.3.2 of~\cite{nagao-tsushima}), the elementary divisors of the Cartan matrix of $G$ are equal to $|C_G(x)|_p$, for $x\in [G_{p'}]$. It follows that the rank of the Cartan matrix modulo $p$, is equal to the number of conjugacy classes of elements of defect~0 of $G$, i.e. the cardinality of $[G_0]$. The following can be viewed as an explicit form of this result:
\begin{enonce}{Proposition}\label{basis}\begin{enumerate}
\item Let $T$ be a simple $kG$-module. Then, in $kR_k(G)$,
$$k\sfc^G(\mathsf{P}_T)=\sum_{x\in[G_0]}\rho\big(\Phi_T(x)\big)\,\gamma_{G,\,x}.$$
\item The elements $\gamma_{G,\,x}$, for $x\in[G_0]$, form a basis of $k\Cart(G)\leq kR_k(G)$.
\end{enumerate}
\vspace{-3ex}
\end{enonce}
\pf Throughout the proof, we simply write $\gamma_x$ instead of $\gamma_{G,\,x}$.\mpn
1. By~\ref{rho}, we have
\begin{align*}
k\sfc^G(\mathsf{P}_T)&=\sum_{S\in\Irr(kG)}\rho(\sfc^G_{T,S})S=\sum_{S\in\Irr(kG)}\sum_{x\in [G_0]}\rho\left(\frac{\Phi_T(x)\Phi_S(x^{-1})}{|C_G(x)|}\right)S\\
&=\sum_{x\in [G_0]}\sum_{S\in\Irr(kG)}\rho\left(\frac{\Phi_T(x)\Phi_S(x^{-1})}{|C_G(x)|}\right)S\\
&=\sum_{x\in [G_0]}\rho\big(\Phi_T(x)\big)\sum_{S\in\Irr(kG)}\rho\left(\frac{\Phi_S(x^{-1})}{|C_G(x)|}\right)S\\
&=\sum_{x\in[G_0]}\rho\big(\Phi_T(x)\big)\,\gamma_x.
\end{align*}
2. We first prove that $\gamma_x$ lies in the image of $k\sfc^G$, for any $x\in G_0$. So let $x\in\nolinebreak G_0$, and $1_x:\langle x\rangle \to \CO$ be the map with value 1 at $x$ and~0 elsewhere. Then $|x|1_x=\sum_{\zeta}\limits\zeta(x^{-1})\zeta$, where $\zeta$ runs through the simple $k\langle x\rangle$-modules, i.e. the group homomorphisms $\langle x\rangle\to \CO^\times$, is an element of $\CO P_k(\langle x\rangle)=\CO R_k(\langle x\rangle)$. Let $v_x=\Ind_{\langle x\rangle}^G(|x|1_x)$. Then $v_x\in \CO P_k(G)$, and its modular character evaluated at $g\in G$ is equal to
\begin{align}\refstepcounter{monequation}
\Phi_{v_x}(g)&=\frac{1}{|x|}\sum_{\substack{h\in G\\g^h\in\langle x\rangle}}\Phi_{|x|1_x}(g^h)\nonumber\\
\label{character}&=\frac{1}{|x|}\sum_{\substack{h\in G\\g^h=x}}|x|=\left\{\begin{array}{ll}|C_G(x)|&\hbox{if}\;g=_Gx\\0&\hbox{otherwise},\end{array}\right.
\end{align}
where $g=_Gx$ means that $g$ is conjugate to $x$ in $G$. Now from Assertion 1, we get that 
\begin{moneq}\label{v viz gamma}k\sfc^G(v_x)=\sum_{y\in[G_0]}\rho\big(\Phi_{v_x}(y)\big)\gamma_y=|C_G(x)|\gamma_x,
\end{moneq}
so $\gamma_x$ is in the image of $k\sfc^G$, since $|C_G(x)|\neq 0$ in $k$.\par
Now by Assertion 1, the elements $\gamma_x$, for $x\in[G_0]$, generate the image $\Cart(G)$ of $k\sfc^G$. They are moreover linearly independent: Suppose indeed that some linear combination $\sum_{x\in [G_0]}\limits\lambda_x\gamma_x$, where $\lambda_x\in k$, is equal to 0. For all $x\in\nolinebreak{[G_0]}$, choose $\widetilde{\lambda}_x\in \CO$ such that $\rho(\widetilde{\lambda}_x)=\lambda_x$. By~(\ref{v viz gamma}), we get an element $\sum_{x\in[G_0]}\limits \widetilde{\lambda}_x\frac{v_x}{|C_G(x)|}$ of $\CO\Proj(kG)$ whose modular character has values in the maximal ideal $J(\CO)$ of $\CO$. But by~(\ref{character}), the value at $g\in G_{p'}$ of this modular character is equal to
$$\sum_{x\in[G_0]}\limits \widetilde{\lambda}_x\frac{\Phi_{v_x}(g)}{|C_G(x)|},$$
which is equal to 0 if $g\notin G_0$, and to $\widetilde{\lambda}_x$ if $g$ is conjugate to $x\in [G_0]$ in $G$. It follows that $\widetilde{\lambda}_x\in J(\CO)$, hence $\lambda_x=\rho(\widetilde{\lambda}_x)=0$. Since $g\in G_{p'}$ was arbitrary, we get that $\lambda_x=0$ for any $x\in [G_0]$, so the elements $\gamma_x$, for $x\in [G_0]$, are linearly independent. This completes the proof of Proposition~\ref{basis}.\endpf
\tpar {\bf The functors $\SS_{L,1,W}$.} In this section, we consider the case where $u=1$, i.e. the case of simple functors $\SS_{L,1,W}$, where $L$ is a $p$-group and $W$ is a simple $\F \Out(L)$ module. The case where $\F$ has characteristic 0 is solved by Corollary~7.4 of~\cite{functorialequivalence}. Then we are left with the case where $\F$ has characteristic~$p$, and we assume that $\F=k$. In this situation, for a finite group $G$, the set $\CP(G,L,1)$ is just the set of pairs $(Q,\delta)$, where $Q$ is a $p$-subgroup of $G$ and $\delta:L\to Q$ is a group isomorphism. Moreover, the set $[\CP(G,L,1)]$ is in one to one correspondence with a set of representatives of conjugacy classes of subgroups $Q$ of $G$ which are isomorphic to $L$ (the bijection mapping $(Q,\delta)$ to $Q$). Then the group $G_{Q,\delta,1}$ is just $N_G(Q)$, while $G_{Q,\delta}=QC_G(Q)$. By Theorem~\ref{the simple functors}, we have that
$$\SS_{L,1,W}(G)\cong  \mathop{\bigoplus}_{(Q,\delta)\in [\CP(G,L,1)]}\limits\mathrm{Tr}_{\un}^{N_G(Q)/QC_G(Q)}\Big(k\Cart\big(C_G(Q)/Z(Q)\big)\otimes W\Big).$$
\begin{enonce}{Notation} Let $G$ be a finite group.
\begin{itemize} 
\item For a subgroup $Q$ of $G$ and an element $z$ of $G$, we set $N_G(Q,z)=N_G(Q)\cap C_G(z)$ and $C_G(Q,z)=C_G(Q)\cap C_G(z)$
\item For a $p$-subgroup $Q$ of $G$, we denote by $\zeta(G,Q)$\index{\ordidx{ZB}$\zeta(G,Q)$} the set of elements~$z$ in $C_G(Q)_{p'}$ for which $Z(Q)$ is a Sylow $p$-subgroup of $C_G(Q,z)$. The group $N_G(Q)$ acts on $\zeta(G,Q)$ by conjugation, and we denote by $[\zeta(G,Q)]$ a set of representatives of orbits under this action.
\item For a finite $p$-group $L$, we denote by $\CZ(G,L)$\index{\ordidx{ZA}$\CZ(G,L)$} the set of pairs $(Q,z)$, where $Q$ is a subgroup of $G$ isomorphic to $L$, and $z$ is an element of $\zeta(G,Q)$. In other words 
$$\CZ(G,L)\!=\!\big\{(Q,z)\mid L\cong Q\leq G,\,z\in C_G(Q)_{p'},\,Z(Q)\in\mathrm{Syl}_p\big(C_G(Q,z)\big) \big\}.$$
The group $G$ acts on $\CZ(G,L)$ by conjugation, and we denote by $[\CZ(G,L)]$ a set of representatives of orbits under this action.
\end{itemize}
\vspace{-2ex}
\end{enonce}
\begin{enonce}{Theorem} \label{simple L,1,W}Let $L$ be a $p$-group and $W$ be a simple $k\Out(L)$-module. Let moreover $G$ be a finite group.
\begin{enumerate}
\item Let $Q$ be a $p$-subgroup of $G$. Then there is an isomorphism $$k\Cart\big(C_G(Q)/Z(Q)\big)\cong \bigoplus_{z\in[\zeta(G,Q)]}\Ind^{N_G(Q)/QC_G(Q)}_{N_G(Q,z)C_G(Q)/QC_G(Q)}k$$
of $kN_G(Q)/QC_G(Q)$-modules.
\item The evaluation of the simple functor $\SS_{L,1,W}$ at $G$ is 
$$\SS_{L,1,W}(G)\cong\bigoplus_{(Q,z)\in[\CZ(G,L)]}\mathrm{Tr}_\un^{N_G(Q,z)/QC_G(Q,z)}(W).$$
\end{enumerate}
\end{enonce}
\vspace{-2ex}
\pf 1. For a subgroup $Q$ of $G$, denote by $x\mapsto \sur{x}$ the projection map $N_G(Q)\to \sur{N}_G(Q)=N_G(Q)/Q$, and set $\sur{C}_G(Q)=QC_G(Q)/Q\cong C_G(Q)/Z(Q)$.\par
By Proposition~\ref{basis}, the vector space $k\Cart\big(\sur{C}_G(Q)\big)$ has a basis consisting of the elements $\gamma_{\sur{C}_G(Q),x}$, for $x\in [\sur{C}_G(Q)_0]$, and the group $N_G(Q)$ permutes these elements. Now if $x\in \sur{C}_G(Q)_{p'}$, there is an element $z\in \big(QC_G(Q)\big)_{p'}$ such that $x=\sur{z}$, and we can moreover assume that $z\in C_G(Q)_{p'}$. Then the centralizer of $\sur{z}$ in $\sur{C}_G(Q)$ is equal to $QC_{QC_G(Q)}(z)/Q=QC_G(Q,z)/Q$. So $\sur{z}\in \sur{C}_G(Q)_0$ if and only if $Q$ is a Sylow $p$-subgroup of $QC_G(Q,z)$, or equivalently, if $Z(Q)$ is a Sylow $p$-subgroup of $C_G(Q,z)$, that is $z\in\zeta(G,Q)$. \par
Moreover, an element $n\in N_G(Q)$ stabilizes $\gamma_{\sur{C}_G(Q),\sur{z}}$ if and only if $\gamma_{\sur{C}_G(Q),\sur{z}}=  \gamma_{\sur{C}_G(Q),\sur{nzn^{-1}}}$, hence by Proposition~\ref{basis}, if there exists $\sur{c}\in \sur{C}_G(Q)$ such that $\sur{nzn^{-1}}=\sur{c}\sur{z}\sur{c}^{-1}$. In other words $\sur{c^{-1}n}\in C_{\sur{N}_G(Q)}(\sur{z})=N_G(Q,z)/Q$, where we set $N_G(Q,z)=N_G(Q)\cap C_G(z)$. So the stabilizer of $\gamma_{\sur{C}_G(Q),\sur{z}}$ in $N_G(Q)$ is equal to $QC_G(Q)N_G(Q,z)=N_G(Q,z)C_G(Q)$. \par
Hence $k\Cart\big(\sur{C}_G(Q)\big)$ is isomorphic to the permutation $N_G(Q)/QC_G(Q)$-module $k\zeta(G,Q)$. The elements $\gamma_{\sur{C}_G(Q),\sur{z}}$, for $z\in[\zeta(G,Q)]$ form a set of representatives of orbits for the action of $N_G(Q)/QC_G(Q)$, and the stabilizer of $\gamma_{\sur{C}_G(Q),\sur{z}}$ is the group $N_G(Q,z)C_G(Q)/QC_G(Q)$. This proves Assertion 1.\spn
2. Now $N_G(Q,z)C_G(Q)/QC_G(Q)\cong N_G(Q,z)/QC_G(Q,z)$, and Assertion 2 follows from Theorem~\ref{the simple functors}, thanks to the general following fact (see Proposition 5.6.10 (ii) in~\cite{benson1}): If $\Gamma'$ is a subgroup of a finite group $\Gamma$, if $M$ is a finite dimensional $k\Gamma$-module and $M'$ is a finite dimensional $k\Gamma'$-module, then
$$\mathrm{Tr}_\un^\Gamma\big((\Ind_{\Gamma'}^\Gamma M')\otimes_k M\big)\cong M'\otimes_k\mathrm{Tr}_\un^{\Gamma'}(\Res_{\Gamma'}^\Gamma M)$$ 
as $k$-vector spaces.  \endpf
\begin{rem}{Remark} The formula in Assertion 2 of Theorem~\ref{simple L,1,W} can be viewed as another instance of similar formulas in Proposition~8.8 of~\cite{thevwebb}, Theorem~2.6 of~\cite{webbglob}, or Theorem 6.1 of~\cite{vanishing}.
\end{rem}
The following corollary deals with the case of Theorem~\ref{simple L,1,W} where $W$ is the trivial module~$k$. First a definition:
\begin{enonce}{Definition} \label{defect isomorphic}Let $G$ be a finite group, and $L$ be a finite $p$-group. An element $z\in G_{p'}$ is said to have {\em defect isomorphic to $L$} if $L$ is isomorphic to a Sylow $p$-subgroup of $C_G(z)$.
\vspace{-2ex}
\end{enonce}
\begin{enonce}{Corollary} \label{simple L,1,k}Let $G$ be a finite group, and $L$ be a finite $p$-group. Then the dimension of $\SS_{L,1,k}(G)$ is equal to the number of conjugacy classes of elements of $G_{p'}$ with defect isomorphic to $L$.
\vspace{-4ex}
\end{enonce}
\pf Indeed, if $(Q,z)\in\CZ(G,L)$ then $\mathrm{Tr}_\un^{N_G(Q,z)/QC_G(Q)}(k)$ is equal to zero if $p$ divides the order of $N_G(Q,z)/QC_G(Q)$, and one dimensional otherwise, that is, if a Sylow $p$-subgroup of $N_G(Q,z)$ is contained in $QC_G(Q,z)$. But since $z\in\zeta(G,Q)$, the group $Q$ is a Sylow $p$-subgroup of $QC_G(Q,z)$. So $\mathrm{Tr}_\un^{N_G(Q,z)/QC_G(Q)}(k)$ is non zero (and then, one dimensional) if and only if $Q$ is a Sylow $p$-subgroup of $N_G(Q,z)=N_{C_G(z)}(Q)$, i.e. if $Q$ is a Sylow $p$-subgroup of $C_G(z)$.\endpf
\vspace{-3ex}
\printindex
\bibliographystyle{abbrv}

\centerline{\rule{5ex}{.1ex}}
\begin{flushleft}
Serge Bouc, CNRS-LAMFA, Universit\'e de Picardie, 33 rue St Leu, 80039, Amiens, France.\\
{\tt serge.bouc@u-picardie.fr}\vspace{1ex}\\
Deniz Y\i lmaz, Department of Mathematics, Bilkent University, 06800 Ankara, Turkey.\\
{\tt d.yilmaz@bilkent.edu.tr}
\end{flushleft}

\end{document}